\documentclass[a4paper,reqno]{amsart}

\usepackage[ngerman,american]{babel}
\usepackage{ucs}
\usepackage[utf8x]{inputenc}
\usepackage{times,euler}
\usepackage{a4wide}
\usepackage{booktabs}
\usepackage{siunitx}
\newif\ifIfAM\IfAMtrue
\ifIfAM\let\qty\SI\fi
\def\isep{\ifIfAM\qquad\else\quad\fi}
\usepackage[capitalize]{cleveref}
\usepackage{ao-abbr}
\usepackage{ao-matrix}
\usepackage{ao-math-fields}
\usepackage{ao-math-symbols}
\usepackage{ao-math-ops}
\usepackage{ao-math-opt}
\usepackage{luh-colors}

\usepackage{caption}
\usepackage{subcaption}
\usepackage{pgfplots}
\usepackage{pgfplotstable}
\pgfplotsset{compat=1.14}
\usetikzlibrary{spy}
\usepgfplotslibrary{colorbrewer}
\usepgfplotslibrary{patchplots}
\usepackage{pgf,tikz}
\usetikzlibrary{arrows}
\usetikzlibrary{arrows.meta}
\usetikzlibrary{positioning}
\usetikzlibrary{fit,calc,patterns,decorations.pathmorphing,decorations.markings}

\usepackage{ifthen}
\usepackage{xifthen}

\usepackage{sl-beam}

\usepackage{fontawesome}

\let\body\Omega
\let\surf\Gamma
\let\map\varphi
\newcommand\id{\textup{id}}

\newcommand\fext{f}
\newcommand\sext{s}

\newcommand\Sym[1]{\mathbb S^{#1}}
\newcommand\TE{E}
\newcommand\TS{\Sigma}

\newcommand\pd{\partial}
\newcommand\pme{\varphi_e}
\newcommand\pms{\varphi_s}
\newcommand\xe{\~e}
\newcommand\xs{\~s}
\newcommand\ds{\~Z}
\newcommand\Phid{\Phi_0}
\newcommand\Phil{\_\Phi}

\newcommand\K{\mathcal K}
\newcommand\M{\mathcal M}
\newcommand\setN{\_\setA}
\renewcommand\cI{c}
\newcommand\cN{c_{\_\setA}}
\newcommand\xiN{\xi_{\_\setA}}
\newcommand\xiA{\xi_\setA}
\newcommand\range[2]{\set{#1,\dots,#2}}

\newtheorem{prop}{Proposition}
\theoremstyle{remark}
\newtheorem{remark}{Remark}
\newtheorem*{remarks}{Remarks}

\DeclareMathOperator\Tr{tr}

\let\mc\multicolumn
\def\tQP{1.49e-08}
\def\E#1#2{1.00e-#1#2}
\def\ok{$\checkmark$}
\def\badobj{objective $>0$}
\def\badcc{$c(q)\tp \xi \ne 0$}

\begin{document}

\title
[Formulating and Heuristic Solving of Contact Problems in Hybrid DDCM]
{Formulating and Heuristic Solving of Contact Problems\\
  in Hybrid Data-Driven Computational Mechanics}

\author[C.G. Gebhardt]{Cristian G. Gebhardt}
\address{Cristian G. Gebhardt\\
  University of Bergen\\
  Geophysical Institute and Bergen Offshore Wind Centre (BOW)\\
  Allégaten 70\\5007 Bergen\\Norway}
\email{cristian.gebhardt@uib.no}
\urladdr{uib.no/en/persons/Cristian.Guillermo.Gebhardt}

\author[S. Lange]{Senta Lange}
\address{Senta Lange\\
  Leibniz Universität Hannover\\Institute of Applied Mathematics\\
  Welfengarten 1\\30167 Hannover\\Germany}
\email{lange@ifam.uni-hannover.de}
\urladdr{ifam.uni-hannover.de/lange}

\author[M.C. Steinbach]{Marc C. Steinbach}
\address{Marc C. Steinbach\\
  Leibniz Universität Hannover\\Institute of Applied Mathematics\\
  Welfengarten 1\\30167 Hannover\\Germany}
\email{mcs@ifam.uni-hannover.de}
\urladdr{ifam.uni-hannover.de/mcs}

\begin{abstract}
  In this work we consider the
  hybrid Data-Driven Computational Mechanics (DDCM) approach,
  in which a smooth constitutive manifold is reconstructed
  to obtain a well-behaved nonlinear optimization problem (NLP)
  rather than the much harder discrete-continous NLP (DCNLP)
  of the direct DDCM approach.
  The key focus is on the addition of geometric inequality constraints
  to the hybrid DDCM formulation.
  Therein, the required constraint force leads to a contact problem
  in the form of a mathematical program with complementarity constraints (MPCC),
  a problem class that is still less complex than the DCNLP.
  For this MPCC we propose a heuristic quick-shot solution approach,
  which can produce verifiable solutions by solving up to four NLPs.
  We perform various numerical experiments
  on three different contact problems of increasing difficulty
  to demonstrate the potential and limitations of this approach.
\end{abstract}


\date\today

\maketitle

\section{Introduction}
\label{sec:intro}
Solving static or dynamic contact problems is of great interest in many
engineering applications, and substantial efforts have been made to study these
problems and to solve them numerically 
\cite{Laursen_Maker:1995,Wriggers:2006,Temizer:2014,Cantin:2021}.
The focus of this article is on a general strategy for solving
problems with inequalities, specifically contact problems,
whose formulation is based on the latest paradigm in continuum mechanics
known as \emph{Data-Driven Computational Mechanics}, or DDCM in brief.
Traditionally, the relation between strains and stresses
of an elastic material is described by some material model
that has to be calibrated from experimental measurements and
is the sole empirical equation in the overall mechanical model.
In data-rich situations
where a large set of measured strain-stress pairs is available,
the DDCM paradigm replaces the empirical strain-stress equation
by the requirement that strain-stress pairs computed in the overall model
are closest to the given data set.
A spatial discretization of such a principle by finite elements then yields
an equality-constrained optimization model
instead of the traditional set of equations.
Because of the discrete data set, that model is a
discrete-continuous nonlinear optimization problem (DCNLP in brief).
In contrast to well-posed MINLP models, however,
the DCNLP has in general no structure with respect to the discrete data points
from which relaxations could be derived.
Therefore, it is typically solved by meta-heuristic approaches
rather than rigorous optimization methods.
Subsequently, we will refer to the above strategy as the
\emph{direct} DDCM approach
\cite{Kirchdoerfer2016,Kirchdoerfer2017}.
An alternative \emph{inverse} DDCM approach
employs the data set provided to reconstruct a traditional material model
that expresses the stresses as explicit functions of the strains
by means of an energy functional
\cite{Ibanez2017,Ibanez2018a,Ibanez2019}.
In order to combine the strengths of the two approaches
while mitigating their weaknesses,
we have proposed a \emph{hybrid} DDCM approach
in \cite{Gebhardt:2020a,Gebhardt:2020b}.
This approach allows for non-traditional implicit material models
given by a smooth \emph{constitutive manifold}
that has to be reconstructed from the data set
in a first (off-line) step. Nevertheless, no special functional structure
is considered, \ie, the existence of an energy function is not assumed.
In the subsequent (on-line) simulation, we require closeness
of the computed strain-stress pairs to that manifold,
thus obtaining a well-behaved NLP model
that is amenable to rigorous derivative-based algorithms.
As the formulation is already established,
we need to connect this with our extension for contact related problems.
The latter one is really the novelty of the present work.
The basic observation is that the hybrid DDCM approach
naturally extends to problems with inequality constraints,
particularly geometric inequalities and hence contact problems.
In the following, we elaborate on a heuristic approach to solving
contact problems in hybrid DDCM formulation, demonstrating why and how
our hybrid approach provides a highly suited universal framework for
static and dynamic simulations in nonlinear elasticity.

The outline of this article is as follows.
\Cref{sec:basics} provides an overview of the underlying mathematical setting
by establishing the terminology and notation for
the continuous and discrete boundary value problem,
the finite element formulations used in this work,
and eventually the DDCM formulations.
Furthermore, the model of the geometrically exact beam is introduced concisely
to be used in the subsequent numerical experiments.
In \cref{sec:contact}, we present the extension of the hybrid DDCM formulation
to contact problems and propose a heuristic NLP-based approach
for solving the resulting MPCC.
We then provide an extensive numerical investigation of our approach
on three different contact problems for the geometrically exact beam model
in \cref{sec:example}.
Finally, in \cref{sec:conclusion},
we close with a summary and some final remarks.

\section{Mathematical Setting}
\label{sec:basics}
\subsection{Continuous boundary value problem}
\label{sec:conti-BVP}

Although we are primarily interested in spatially discretized models,
we need to fix the underlying continuous setting.
We brief\/ly describe a simplified static model problem
(omitting most technical details and
assuming sufficient regularity throughout)
and give some comments on typical generalizations afterward.
The complete theory can be found, \eg, in
\cite{Ciarlet:1988,Ciarlet:2005}.

Given a \emph{domain} $\body \subset \R^3$
(a bounded and connected open subset
with a Lipschitz-continuous boundary $\surf = \boundary \body$),
consider an elastic body that occupies the closure $\_\body$
(its \emph{reference configuration}) and that is subjected to
\emph{body forces} $\fext\: \body \to \R^3$ and
\emph{surface forces} $\sext\: \surf \to \R^3$.
The induced \emph{deformation} $\map:\_\body \to \R^3$
is an orientation-preserving diffeomorphism
($\det \map'(x) > 0$ for every $x \in \_\body$)
whose image $\map(\_\body) \subset \R^3$ is called
the body's \emph{deformed configuration}.
An equivalent description is given by the
\emph{displacement field} $u\:\_\body \to \R^3$
where $\map = \id + u$.
We seek the unknown displacement field $u$
in a suitable function space $V$
as solution of the weak formulation
of the \emph{equations of equilibrium},
called the \emph{principle of virtual work} in this context.
The principle of virtual work states that there exists
a symmetric matrix field $\TS\: \_\body \to \Sym3 \subset \R^{3 \x 3}$,
known as the \emph{second Piola--Kirchhoff stress tensor field}, such that
\begin{equation}
  \label{eq:virt-work}
  \int_\body (I + \nabla u(x)) \TS(x) : \nabla v(x) \, dx
  =
  \int_\body \sprod{\fext(x)}{v(x)} \, dx
  +
  \int_\surf \sprod{\sext(x)}{v(x)} \, dS
  \quad \forall v \in V
  .
\end{equation}
Here we use the common notation $\nabla u(x) = u'(x) \in \R^{3 \x 3}$
and $A : B = \Tr(A\tp B)$ for $A, B \in \R^{3 \x 3}$.
The first integrand can be rewritten as
$(I + \nabla u) \TS : \nabla v = \TS : \TE'(u) v$
where $\TE'(u) v$ is the (Gâteaux) derivative of the symmetric
\emph{Green--St.~Venant strain tensor field}
$\TE(u)\: \_\body \to \Sym3$ associated with $u \in V$,
\begin{equation*}
  \TE(u)
  \define
  \frac12 (\nabla u\tp + \nabla u + \nabla u\tp \nabla u)
  =
  \frac12 (\nabla\map\tp \nabla\map - I)
  .
\end{equation*}
In the special case of a \emph{hyperelastic} material, there exists a
\emph{stored energy function} $W\: \_\Omega \x \Sym3 \to \R$
such that $\TS(x) = \pd_{\TE} W(x, \TE(u)(x))$.
This gives the total derivative
$\TS(x) : \TE'(u)(x) v = \pd_u W(x, \TE(u)(x)) v$,
and the principle of virtual work
takes the simple variational form
\begin{equation}
  \label{eq:virt-work-var}
  J'(u) v = 0 \quad \forall v \in V,
\end{equation}
which characterizes $u \in V$ as a stationary point
of the \emph{energy} functional
\begin{equation}
  \label{eq:primal-statement}
  J(u)
  \define
  \int_\body W(x, \TE(u)(x)) \, dx - \left(
  \int_\body \sprod{\fext(x)}{v(x)} \, dx +
  \int_\surf \sprod{\sext(x)}{v(x)} \, dS \right)
  .
\end{equation}
The hyperelastic case is in fact the ``traditional'' case
mentioned in the introduction.
In general, elastic materials are defined by the existence
of a \emph{response function} $R\: \_\Omega \x \Sym3 \to \Sym3$
such that the \emph{constitutive equation} holds:
\begin{equation}
  \label{eq:response}
  \TS(x) = R(x, \TE(u)(x)).
\end{equation}
This explains the existence of $\TS$ in \eqref{eq:virt-work}.
The complete boundary value problem then consists of
\eqref{eq:virt-work} and \eqref{eq:response}
or of \eqref{eq:virt-work-var} and $\TS = \pd_{\TE} W$.
In the discrete setting, we will address an even more general material model
where $\TS$ and $\TE(u)$ are linked
by an implicit constitutive equation.

\begin{remarks}
  The static model \eqref{eq:virt-work} is simplified insofar as
  the body is generally fixed at a subset $\surf_0$ of $\surf$
  while surface forces are only applied at the complement $\surf_1$, \ie,
  $\surf = \surf_0 \cup \surf_1$ with
  $u = 0$ on $\surf_0$ and $\sext\:\surf_1 \to \R^3$.
  Without the fixing, unbalanced forces will physically accelerate the body
  and the static model will not admit a solution.
  Naturally, everything becomes time-dependent in a dynamic model
  and the purely spatial elliptic PDE behind \eqref{eq:virt-work}
  is replaced by a parabolic or hyperbolic PDE
  with initial conditions in time.
  However, the overall model structure does not change fundamentally,
  and we will restrict ourselves to the static case in the following.
  Another major generalization allows hybrid mechanical systems that
  consist of several elastic and possibly rigid bodies linked by joints.
  While systems of this type are rarely considered in a continuous setting,
  they are of course omnipresent in real-life applications,
  and we will in fact address them in the discrete setting.
\end{remarks}

In certain situations, it is more convenient to replace
the (primal) statement (\ref{eq:primal-statement})
by a more general (primal-dual) statement given by
\begin{align}
  \label{eq:primal-dual-statement}
  L(u,\TE,\TS)
  \define
  \int_\body W(x, \TE(x)) \, dx
  &- \int_\body \sprod{\TS(x)}{\TE(x)-\TE(u)(x)} \, dx \\
  \notag
  &- \left(
    \int_\body \sprod{\fext(x)}{v(x)} \, dx +
    \int_\surf \sprod{\sext(x)}{v(x)} \, dS \right)
  ,
\end{align}
which upon a classical stationarity argument, delivers automatically:
the weak form of the equation of equilibrium;
the compatibility condition relating displacement and strain fields,
\ie, $\TE - \TE(u) = 0 $; and,
the constitutive equation relating strain and stress fields,
\ie, $\TS - \pd_\TE W(\fcdot, \TE) = 0$.
This is the type of formulation that we are going to use in the following.

\subsection{Discrete boundary value problem}
\label{sec:discr-BVP}

Consider now a finite element discretization of an elastic body
described by variable vectors $(q, e, s) \in \R^{n_q + n_e + n_s}$.
The \emph{generalized coordinates}~$q$ capture the body's deformation,
expressed, for instance, in terms of
positions in $\R^3$ and orientations in $SO(3)$ or in $S^2$
of all finite element nodes.
Special coordinates may be chosen, \eg, when modeling
essentially one-dimensional bodies (\emph{beams}) and
essentially two-dimensional bodies (\emph{shells})
\cite{Ciarlet:1997,Ciarlet:2000}.
The strain and stress vectors $e$ and~$s$ (with $n_s = n_e$)
collect the strain and stress tensors $\TE$ and $\TS$
at all finite element nodes
in so-called \emph{Voigt notation}
which stores just the relevant entries of each symmetric tensor
in a specific order \cite{Helnwein:2001}.
Formal details of $e$, $s$ are given in the appendix,
details of $q$ are irrelevant here.

In the general case of hybrid mechanical systems
we have one set of all variable types for each elastic body
and one additional pair of position and orientation coordinates
for each rigid body.
Moreover, the coupling of bodies by \emph{joints}
is modeled with \emph{kinematic constraints}.
Kinematic constraints will also appear in the single-body case
when \emph{redundant coordinates} are used.
For instance, an orientation $D$ in the 3-manifold $SO(3)$
is often represented by the nine components
of three \emph{directors} $d_i$ (the columns of the matrix $D$)
together with six orthonormality constraints
$\sprod{d_i}{d_j} = \delta_{ij}$, or $D\tp D = I$ in brief.

The model equations with $n$ degrees of freedom
plus $m$ redundant coordinates ($n_q = n + m$) read
\begin{align}
  \label{eq:kinematic}
  h(q) &= 0 \in \R^m && \text{kinematic constraint}, \\
  \label{eq:compatibility}
  e(q) - e &= 0 \in \R^{n_e} && \text{compatibility}, \\
  \label{eq:equilibrium}
  N(q)\tp (B(q)\tp s - f) &= 0 \in \R^n && \text{equilibrium}, \\
  \label{eq:constitutive}
  g(e, s) &= 0 \in \R^{n_e} && \text{constitutive equation},
\end{align}
where $B(q) \define e'(q) \in \R^{n_e \x n_q}$,
and $N(q) \in \R^{n_q \x n}$ is a nullspace basis
of $H(q) \define h'(q) \in \R^{m \x n_q}$.
Thus we require $\rank N(q) = n$ and hence $\rank H(q) = m$.
Moreover, we assume that $\rank g'(e, s) = n_e$.
The structure of $g$ depends on the structure of $e, s$
and is again described in the appendix.
\begin{remark}
  If $h(q) = 0$ admits a global representation
  $q(x)$ with $x \in \R^{n_q - m}$,
  the reduced equilibrium with $q(x)$ substituted
  is called a \emph{minimal irreducible form} of the static problem.
\end{remark}

Note that (with a slight abuse of notation)
we have introduced a separate strain variable $e$
in addition to the functional dependence $e(q)$
that represents $\TE(u)$.
Note also that the force equilibrium derived from \eqref{eq:virt-work}
lives in the tangent space of the $n$-dimensional
submanifold of $\R^{n_q}$ that the kinematic constraint $h(q) = 0$ defines.
Here $B(q)\tp s = e'(q)\tp s$ corresponds to $\TS : E'(u)$
and $f \in \R^{n_q}$ represents body and surface forces.
Note finally that the implicit constitutive equation (the material model)
defines a \emph{constitutive manifold}
to which the pairs $(e, s)$ have to belong.

\subsection{Standard and general finite element formulations}
\label{sec:DDCM-FEM}

With additional variables $(\xe, \xs) \in \R^{2 n_e}$
that will be relevant later in the DDCM setting,
we now define mappings
\begin{align*}
  \Phid\: &\R^{n_q + 2 n_e} \to \R^{n_q + n_e}, &
  \Phid(q, e, s) &\define
  \col(c){h(q) \\ e(q) - e \\ N(q)\tp (B(q)\tp s - f)}, \\
  \Phil\: &\R^{n_q + 4 n_e} \to \R^{n_q + 2 n_e}, &
  \Phil(q, e, s, \xe, \xs) &\define
  \col(c){\Phid(q, e, s) \\ g(\xe, \xs)}, \\
  \Phi\: &\R^{n_q + 2 n_e} \to \R^{n_q + 2 n_e}, &
  \Phi(q, e, s) &\define \Phil(q, e, s, e, s).
\end{align*}
The \emph{general finite element model}
\eqref{eq:kinematic}--\eqref{eq:constitutive}
then takes the simple form
\begin{equation*}
  \Phi(q, e, s) = 0,
  \tag{gFEM}
\end{equation*}
and it specializes to the \emph{standard finite element model}
(applicable to hyperelastic materials only)
where $s$ is an explicit function of $e$ via $g(e, s) = \psi'(e) - s$
with $\psi(e)$ representing the stored energy $W(x, E)$,
\begin{equation*}
  \Phi(q, e, \psi'(e)) = 0.
  \tag{sFEM}
\end{equation*}

For numerical stability and uniqueness of solutions
we need full rank of the Jacobians
$\Phi', \Phil'$ at feasible points,
where $\Phi(q, e, s) = 0$ or $\Phil(q, e, s, \xe, \xs) = 0$ holds, respectively.
The associated Jacobians read
\begin{align*}
  \Phi'(q, e, s)
  &=
  \mat[5]{
    H(q) \\
    B(q) & -I \\
    D(q, s) & 0 & N(q)\tp B(q)\tp \\
    0 & G_1(e, s) & G_2(e, s)
  }, \\
  \Phil'(q, e, s, \xe, \xs)
  &=
  \mat[5]{
    H(q) \\
    B(q) & -I \\
    D(q, s) & 0 & N(q)\tp B(q)\tp \\
    0 & 0 & 0 & G_1(\xe, \xs) & G_2(\xe, \xs)
  },
\end{align*}
where $G_j(e, s) \define \pd_j g(e, s)$ and
\begin{equation*}
  D(q, s) \define N(q)\tp \pd_q (B(q)\tp s) + W(q, B(q)\tp s - f)
  \qtextq{with}
  W(q, a) \define \pd_q (N(q)\tp a)
  .
\end{equation*}
We note that $W(q, a) = 0$ if $H(q)$ and hence $N(q)$ is constant,
and that $\pd_q (B(q)\tp s)$ is symmetric.

\begin{prop}
  \label{prop:rank}
  If $\mat[2]{D(q, s) N(q) & N(q)\tp B(q)\tp}$ has full rank $n$,
  then the Jacobian $\Phil'(q, e, s, \xe, \xs)$ has full rank $n_q + 2 n_e$,
  \ie, $\dim \ker \Phil'(q, e, s, \xe, \xs) = 2 n_e$,
  and the square Jacobian $\Phi'(q, e, s)$ is invertible.
\end{prop}

We omit the proof, which is a straightforward computation.

\subsection{Simple examples on the real line}

To illustrate the above ideas, we present next some very simple examples
 comprising a single two-node element.
Consider $q = (x_1, x_2) \in \R^2$
where $x_1$ is the displacement coordinate of the left end of a horizontal rod
and $x_2$ is the displacement of its right end.
In addition, if we require the left end to be fixed at the origin ($x_1 = 0$),
we have the linear kinematic constraint
\begin{align*}
	h(q) &= H q = 0, &
	H &= \mat[2]{1 & 0}, &
	N &= \mat[1]{0 \\ 1}.
\end{align*}
Now, we are required to choose a potential energy.
The simplest convex function that can be adopted for the potential energy is
\begin{equation*}
	\psi(e) = \frac12 c e^2,
\end{equation*}
where $c$ is a positive constant
that depends on the material and geometrical properties.
Such a potential energy automatically imposes
\begin{equation*}
	s = \psi'(e) = c e \qtextq{and} k = \psi''(e) = c,
\end{equation*}
where $k$ is usually known as the stiffness.
This can also be understood as a form of Hooke's law.
Lastly, we consider the vertical external force $f = (0, f_2)$.

\paragraph{Example 1: the simplest linear case}

In this first example, we consider a deformation measure
given by the linear expression
\begin{equation*}
	e(q) = \frac{x_2 - x_1}{l} - 1,
\end{equation*}
where the constant $l$ corresponds to the value of $x_2-x_1$ for a vanishing $f$,
\ie, the natural state of the system.
The deformation measure proposed leads to
\begin{equation*}
	E = \frac1l \mat[2]{-1 & 1}.
\end{equation*}

Now, by computing $s(q) = c e(q)$ and inserting it
into the the reduced balance equation,
we obtain the following system of linear equations:
\begin{align*}
	N\tp (E\tp s(q) - f) &= 0\: &
	\frac{c}{l} \left( \frac{x_2 - x_1}{l} - 1 \right) - f_2 = 0,\\
	h(q) &= 0\: & x_1 = 0,
\end{align*}
which can be further simplified to its minimal irreducible form
\begin{equation*}
	x_2 - l \left( 1 + \frac{l f_2}{c} \right) = 0.\\
\end{equation*}
From this equation we can trivially conclude
that if $f_2$ is positive (traction force),
$x_2-x_1$ becomes larger than $l$ and thus,
the system is undergoing stretching.
Conversely, if $f_2$ is negative (compression force),
$x_2-x_1$ becomes smaller than $l$ and thus,
the systems is undergoing shortening.

The two matrices considered in \cref{prop:rank}
(with variables omitted) have obviously full rank:
\begin{align*}
	\mat[2]{D N & N\tp B\tp} &= \mat[2]{0 & 1/l},
	&
	\Phil' &= \mat[6]{
		1 & 0 \\
		-1/l & 1/l & -1 \\
		0 & 0 & 0 & 1/l \\
		0 & 0 & 0 & 0 & 0 & -1}.
\end{align*}

\begin{remark}
	The term $(x_2 - x_1) / l$ present in the definition of $e(q)$
	is not good enough to measure the deformation of the structure,
	although being invariant under translations.
	For a proper physical measure, $E$ has to be a difference
	of two Riemannian metrics.
	Such a problem is easily circumvented through a term of the form
	$\frac12 (x_2 - x_1)^2 / l^2$.
\end{remark}

\paragraph{Example 2: the simplest nonlinear case}

In this second example, we consider a deformation measure
given by the quadratic expression
\begin{equation*}
	e(q) = \frac12 \left( \left( \frac{x_2 - x_1}{l} \right)^2 - 1 \right)
\end{equation*}
that leads to
\begin{equation*}
	E(q) = \frac1l \mat[2]{-x_1 & x_2}.
\end{equation*}
Now, by computing $s(q)$ and inserting it
into the the reduced balance equation,
we obtain the following system of nonlinear equations:
\begin{align*}
	N\tp (E(q)\tp s(q) - f) &= 0\: &
	\frac{c x_2}{2 l} \left( \left(
	\frac{x_2 - x_1}{l} \right)^2 - 1 \right) - f_2 = 0, \\
	h(q) &= 0\: & x_1 = 0,
\end{align*}
which can be further simplified to its minimal irreducible form
\begin{equation*}
	x_2^3 + l^2 x_2^2 + \frac{2 l^3 f_2}{c} = 0.
\end{equation*}
To discuss the meaning of the roots, let us consider
$c = 1$, $l = 1$ and $f_2 = 20$.
Then the previous equation becomes
\begin{equation*}
	x_2^3 + x_2^2 + 40 = 0,
\end{equation*}
with roots $(x_2)_1 \approx +3.3225$, $(x_2)_2 \approx -1.6613+3.0432i$,
and $(x_2)_3 = (\_x_2)_2 \approx -1.6613+3.0432i$.
While the first root is real and positive (which indicates stretching),
the last two are complex conjugate with negative real part
(which indicates shortening and negative length).
As negative length is impossible, the only physically meaningful solution
is then provided by the first purely real root.

Again, the two matrices considered in \cref{prop:rank} have full rank:
\begin{align*}
	\mat[2]{D N & N\tp B\tp} &= \mat[2]{s/l & x_2/l},
	&
	\Phil' &= \mat[6]{
		1 & 0 \\
		-x_1/l & x_2/l & -1 \\
		0 & s/l & 0 & 1/l \\
		0 & 0 & 0 & 0 & 0 & -1}.
\end{align*}

\subsection{Specialization to the geometrically exact beam model}

For this model, the configuration manifold $Q$ is $\R^3 \x SO(3)$.
However and due to several analytical and numerical advantages,
it is more convenient to describe the system's configuration
in a larger manifold $P$ defined as $\R^3 \x \R^{3 \x 3}$.
$P$ is such that $P \supset Q$,
which of course has the very favorable vector space structure
and must be supplemented with certain constraints
to restrict the motion in $P$ to $Q$.
The corresponding configuration map $q: [0, L] \to P$
essentially describes the position of the midline and the orientation
(through three orthonormal directors)
of any cross-section spanned by the coordinates
$(\theta^1, \theta^2) \in \setA(\theta^3)$
located at any length coordinate $\theta^3 \in [0, L]$,
which is given by
\begin{equation*}
  q(\theta^3)
  =
  (x_0(\theta^3), d_1(\theta^3), d_2(\theta^3), d_3(\theta^3))
  \in P \cong \R^{12}.
\end{equation*}
The position map that provides the position
in the ambient space of any point belonging to the beam can be expressed as
\begin{equation*}
  x(\theta)
  =
  x_0(\theta^3) + \theta^1 d_1(\theta^3) + \theta^2 d_2(\theta^3)
  \in \R^3,
\end{equation*}
and the corresponding constraint map is
\begin{equation*}
  h(x_0, d_1, d_2, d_3)
  =
  (\set{\norm{d_i}_2^2 - 1}_{i=1}^3,
  \sprod{d_1}{d_2}, \sprod{d_2}{d_3}, \sprod{d_3}{d_1}) \in \R^6.
\end{equation*}
The model-specific deformation measures are the axial and shear strains
\begin{equation*}
  \Gamma^i \define \sprod{\pd_{\theta^3} x_0}{d_i} - \Gamma^i_0,
\end{equation*}
and the bending and torsional strains
\begin{equation*}
  K^i \define \frac12 \eps^i_{jk} \bigl(
  \sprod{\pd_{\theta^3} d_j}{d_k} -
  \sprod{\pd_{\theta^3} d_k}{d_j} \bigr) - K^i_0.
\end{equation*}

To further describe the model in the continuous setting,
we need to consider the infinite-dimensional space $V \define H^1(0, L)^{12}$
for the configuration map $q$.
Note that, in the present context,
only generalized weak solutions are of our particular interest.
Exact analytical solutions can be found for few simplyfied cases,
but in general one needs to compute discrete solutions.
To this end, we dicretize the beam into $\card{\mathcal N_e}$ intervals
$I_i = [L_{i-1}, L_i]$ such that $\bigcup_{i \in \mathcal N_e} I_i = [0, L]$,
and introduce the discrete space
\begin{equation*}
  V_h \define
  \defset{v_h \in C^0(0, L)^{12} \cap V}
  {v_h \in \Poly^1(I_i) \ \forall i \in \mathcal N_e}.
\end{equation*}
Such a space is where the problem is to be approximated and thus,
by contruction $V_h \subset V$.
Particularly for the geometrically exact beam, we have
\begin{equation*}
  q(\theta^3) \approx q_h(\theta^3)
  =
  \sum_{i \in \mathcal N_e} \Phi_i(\theta^3) \^q_i,
\end{equation*}
where $\Phi_i(\theta^3)$ is a matrix that contains all shape functions
and $\^q_i$ is the set of nodal degrees of freedom.
Upon removal of the rigid body motions
and by considering a standard linear material law
(which is completely omitted in this work),
the optimal convergence rate
\begin{equation*}
  \norm[L^2(0,L)^{12}]{q - q_h} \le c h^2
\end{equation*}
can be achieved, for $h$ being the maximal interval length
among all intervals $I_i$.

This particular discrete model for the geometrically exact beam
satisfies all the conditions requiered by proposition \ref{prop:rank},
see \cite{Gebhardt:2020a,Gebhardt:2020b},
and proves to be sufficiently complex to show
the favorable features of our hybrid DDCM approach.

\subsection{DDCM formulations}

In the DDCM context,
the constitutive equation \eqref{eq:constitutive}
is replaced by a data set $\ds$ of measured
strain-stress pairs $(\xe, \xs)$
where $\ds$ has a product structure
that is described in the appendix.
Then we choose proximity measures $\pme, \pms$
that satisfy $\pme(0) = \pms(0) = 0$,
and the \emph{direct} DDCM approach replacing (gFEM) consists in solving
the \emph{discrete-continuous nonlinear optimization problem}
\begin{equation*}
  \min_{(\xe, \xs) \in \ds, (q, e, s)} \
  \pme(e - \xe) + \pms(s - \xs)
  \qstq \Phid(q, e, s) = 0.
  \tag{DCNLP}
\end{equation*}
Our \emph{hybrid} DDCM approach requires that a constitutive manifold
(defined by $g(\xe, \xs) = 0$)
has been reconstructed from $\ds$ in a first step,
and as replacement of (gFEM) we study the
\emph{approximate nonlinear optimization problem}
\begin{equation*}\label{eq:aNLP}
  \min_{q, e, s, \xe, \xs} \
  \pme(e - \xe) + \pms(s - \xs)
  \qstq \Phil(q, e, s, \xe, \xs) = 0.
  \tag{aNLP}
\end{equation*}
Standard choices of the proximity measures are energy norms induced
by a symmetric positive definite weight matrix $C \in \Sym{n_e}$:
\begin{align*}
  \pme(e - \xe) &= \frac12 \norm[C]{e - \xe}^2, &
  \pms(s - \xs) &= \frac12 \norm[\Inv C]{s - \xs}^2.
\end{align*}

Our focus here is on the nonlinear \emph{optimization problem} (aNLP)
as a replacement of the nonlinear \emph{equation system} (gFEM).
The optimization context permits a direct extension of (gFEM)
to more general problem classes,
in particular to problems that involve inequality constraints.
Below we will address contact problems as a specific case.
On the other hand, the optimization context provides flexibility
in the numerical treatment.
In the following we brief\/ly discuss an important aspect
of this algorithmic flexibility.

Essentially (aNLP) is a kind of penalty formulation of (gFEM)
and thus a pure \emph{feasibility problem}
with $n_q + 2 n_e$ equality conditions
for the same number of variables $(q, e, s)$:
for any solution $(q, e, s)$ of (gFEM),
setting $(\xe, \xs) \define (e, s)$
clearly gives a solution of (aNLP)
with objective value zero and vice versa.
In contrast, solutions of (aNLP) with nonzero objective value
are not solutions of (gFEM):
the computed strain-stress pair $(e, s)$
does not match the closest pair $(\xe, \xs)$
on the constitutive manifold.
Thus, converting the nonlinear equation $g(e,s) = 0$ to the minimization of
$\pme(e - \xe) + \pms(s - \xs)$ naturally permits robust problem formulations
for numerically difficult cases. The simplest robustification just replaces
the $\ell_2$ objective terms by corresponding $\ell_1$ terms,
\begin{align*}
  \pme(e - \xe) &= \norm{C^{1/2} (e - \xe)}_1, &
  \pms(s - \xs) &= \norm{C^{-1/2} (s - \xs)}_1.
\end{align*}
With $z_e \define e - \xe$ and $z_s \define s - \xs$,
the resulting problem is still in the form of (aNLP) and can be written
\begin{equation*}
  \min_{z_e, z_s, q, e, s} \ \norm{C^{1/2} z_e}_1 + \norm{C^{-1/2} z_s}_1
  \qstq \Phid(q, e, s) = 0, \quad g(e - z_e, s - z_s) = 0.
\end{equation*}
This formulation has a clear physical interpretation:
$z_e$ and $z_s$ can be understood as the mininal perturbations
of the strain and stress fields $e$ and $s$
required to satisfy the constitutive equation for any solution $(q, e, s)$
of the balance, compatibility, and kinematic constraints.
An alternative robustification penalizes the residual of $g$
instead of the arguments,
\begin{equation*}
  \min_{z, q, e, s} \ \norm{W_z z}_1
  \qstq \Phid(q, e, s) = 0, \quad g(e, s) = z,
\end{equation*}
where $W_z \in \Sym{n_e}$ denotes a positive definite diagonal weight matrix.
This formulation has a clear physical interpretation as well:
here $z$ can be understood as the minimal relaxation
of the constitutive equation for any solution $(q, e, s)$
of the balance, compatibility, and kinematic constraints.
Finally we might consider the fully robust formulation
\begin{equation*}
  \min_{y, z, q, e, s} \ \norm{W_y y}_1 + \norm{W_z z}_1
  \qstq \Phid(q, e, s) = y, \quad g(e, s) = z
\end{equation*}
with a second residual $y$ and associated weight matrix $W_y$.
Here the interpretation is that $y$ and $z$ are the minimal relaxations
required to make the complete set of constraints feasible
for any triple $(q, e, s)$.
Of course, further formulations with different combinations
of exactly satisfied constraints and relaxed constraints
can be considered when needed.
Finally we notice that all these nonsmooth robust optimization problems
admit smooth NLP reformulations with nonnegative slack variables.

\section{Contact Problems in Hybrid DDCM Formulation}
\label{sec:contact}
\subsection{Formulation of the contact problem}

In this section we consider additional geometric inequality constraints
in (aNLP), that is,
constraints on the configuration variable $q$: $c(q) \ge 0$ in $\R^{n_c}$.
If one of these constraints becomes active, $c_i(q) = 0$,
a corresponding constraint force $\xi_i$ arises,
and hence the balance equation needs to be extended as
\begin{align*}
  N(q)\tp (B(q)\tp s  - f - C(q)\tp \xi).
\end{align*}
Here $C(q) \define c'(q)$ in $\R^{n_c \x n_q}$ is the constraints matrix,
and $\xi \ge 0$ in $\R^{n_c}$ is the vector of all constraint forces
(mathematically a vector of Lagrange multipliers)
which satisfy the additional complementarity condition $c(q)\tp \xi = 0$,
i.e., $c_i(q) \xi_i = 0$ for every $i \in \setI \define \set{1, \dots, n_c}$.
Then (aNLP) generalizes to the following contact problem in form of a
\emph{mathematical program with complementarity constraints (MPCC)}
wherein the (physical) Lagrange multiplier $\xi$ appears as a primal variable:
\begin{align*}
  \label{eq:MPCC}
  \min_{q, e, s, \xe, \xs, \xi} \quad
  & \frac12 \norm[C]{e - \xe}^2 + \frac12 \norm[\Inv C]{s - \xs}^2 \\
  \stq
  & N(q)\tp (B(q)\tp s - f - C (q)\tp \xi) = 0, \\
  & h(q) = 0, \\
  & e(q) - e = 0, \tag{MPCC} \\
  & g(\xe, \xs) = 0, \\
  & c(q) \ge 0, \\
  & \xi \ge 0, \\
  & c(q)\tp \xi = 0.
\end{align*}
It is well-known that the problem class (MPCC) is less regular than (NLP)
because of the complementarity condition $c(q)\tp \xi = 0$:
together with this condition, the entire set of constraints cannot satisfy
a strong \emph{constraint qualification (CQ)} at any feasible point, cf.\
\cite{Luo_et_al:1996,Scheel_Scholtes:2000,Hoheisel_et_al:2013}.
As a consequence, \eqref{eq:MPCC} is generally more difficult to solve
than an NLP.
A typical solution approach converts the MPCC to a parameterized NLP
by regularization and solves a sequence of NLPs
where the regularization parameter is driven to zero.
The regularization can be based on relaxation, penalization, or smoothing,
and warm starting techniques can be used to speed up the solution process, cf.\
\cite{Fukushima_Pang:1999,Jiang_Ralph:2000,Guo_Lin_Je:2015,Wang_Biegler:2023}.

\subsection{A quick shot solution approach}

Here we propose a ``quick shot'' approach using up to four NLPs that can produce
verifiable solutions with moderate effort in fortunate cases.
Even if that first attempt does not succeed,
it is useful because it provides warm start information
for the subsequent sequential MPCC solution procedure.
Our quick shot approach drops complementarity in (MPCC),
fixes $\xi$ or selected parts,
and modifies selected constraints, as follows:
\begin{enumerate}
  \setcounter{enumi}{-1}
\item Just extend (aNLP) with the constraint $c(q) \ge 0$,
  which gives (MPCC) with $\xi$ fixed at zero, and solve.
  If the solution $x^0 = (q^0, e^0, s^0, \xe^0, \xs^0)$ satisfies $c(q^0) > 0$,
  then it is a valid solution of (MPCC) together with $\xi^0 = 0$:
  the underlying assumption that no contact occurs was correct.
\item Otherwise determine the active set (or contact set)
  $\setA \define \defset{i \in \setI}{c_i(q^0) = 0}$
  and denote by $\cA$ and $\xiA$ the associated vectors of
  constraint functions and constraint forces.
  Solve the following NLP using $x^0$ for a warm start and
  initializing $\xiA = 1$ and $\xiN = 0$
  where $\setN \define \setI \setminus \setA$:
  \begin{align*}
    \label{eq:NLP1}
    \min_{q, e, s, \xe, \xs,  \xi} \quad
    & \frac12 \norm[C]{e - \xe}^2 + \frac12 \norm[\Inv C]{s - \xs}^2 \\
    \stq
    & N(q)\tp (B(q)\tp s - f - C(q)\tp \xi) = 0, \\[-2pt]	
    & h(q) = 0, \\[-2pt]
    & e(q) - e = 0, \tag{NLP1} \\[-2pt]
    & g(\xe, \xs) = 0, \\
    & \cA(q) = 0, \cN(q) \ge 0, \\[-2pt]
    & \xiA \ge 0, \ \xiN = 0.
  \end{align*}
  If a solution $(x^1, \xi^1)$ of (NLP1) is found,
  then it is also a solution of (MPCC):
  the underlying assumption that $\setA$ is exactly the contact set was correct.
\item Otherwise relax \eqref{eq:NLP1} by replacing the equality $\cA(q) = 0$
  with the inequality $\cA(q) \ge 0$ to obtain (NLP2),
  and solve using $(x^1, \xi^1)$ for a warm start.
  If the solution $(x^2, \xi^2)$ satisfies the missing
  complementarity condition on $\setA$, $\cA(q^2)\tp \xiA^2 = 0$,
  then it is a solution of (MPCC):
  the underlying assumption that no contact occurs outside $\setA$ was correct.
\item Otherwise relax (NLP2) further by replacing the equality $\xiN = 0$
  with the inequality $\xiN \ge 0$ to obtain (NLP3),
  and solve using $(x^2, \xi^2)$ for a warm start.
  If the solution $(x^3, \xi^3)$ satisfies the missing complementarity
  condition $c(q^3)\tp \xi^3 = 0$, then it is a solution of (MPCC).
\end{enumerate}

\begin{figure}
  \newcommand\drop{\color{LUH-red}}
  \newcommand\impl{\color{LUH-blue}}
  \newcommand\NLPw{6.4em}
  \newcommand\NLPd{2.5em}
  \newcommand\NLP[5]{%
    \node[text width=\NLPw,draw,rectangle,rounded corners=3pt,fill=gray!30]
    (#1) [#3,anchor=north west] {\begin{minipage}{\textwidth}#4\end{minipage}};
    \node[text width=\NLPw] (#2) [below=-1ex of #1.south,anchor=north]
    {\begin{minipage}{\textwidth}\begin{align*}#5\end{align*}\end{minipage}};
  }
  \centering
  \begin{tikzpicture}[->,shorten >=1pt,thin]
    \NLP15{}
    {(aNLP)\\assume:\\[-2pt]no contact\\[-2pt]}
    {c(q) &\ge 0 \\ \\ \impl\xi &\impl= 0 \\ \\ \impl c(q)\tp \xi &\impl= 0 \\}
    \NLP26{right=\NLPd of 1.north east}
    {(NLP1)\\assume:\\[-2pt]contact set = $\setA$\\[-2pt]}
    {\cA(q) &= 0 \\ \cN(q) &\ge 0 \\ \xiA &\ge 0 \\ \impl\xiN &\impl= 0 \\
      \impl\cA(q)\tp \xiA &\impl= 0 \\ \impl\cN(q)\tp \xiN &\impl= 0}
    \NLP37{right=\NLPd of 2.north east}
    {(NLP2)\\assume:\\[-2pt]no contact\\[-2pt]outside $\setA$}
    {\cA(q) &\ge 0 \\ \cN(q) &\ge 0 \\ \xiA &\ge 0 \\ \impl\xiN &\impl= 0 \\
      \drop\cA(q)\tp \xiA &\drop= 0 \\ \impl\cN(q)\tp \xiN &\impl= 0}
    \NLP48{right=\NLPd of 3.north east}
    {(NLP3)\\only use $\setA$\\[-2pt]for warm start\\[-2pt]}
    {c(q) &\ge 0 \\ \\ \xi &\ge 0  \\ \\ \drop c(q)\tp \xi &\drop= 0 \\}
    \draw[-Latex] (1.east) -- (2.west);
    \draw[-Latex] (2.east) -- (3.west);
    \draw[-Latex] (3.east) -- (4.west);
  \end{tikzpicture}
  \caption{The four quick shot NLPs with relevant constraints:
    black explicit, blue implicit, red dropped.
    The splitting of $\setI$ into $\setA$ and $\setN$
    occurs only in (NLP1) and (NLP2).}
  \label{fig:quick-shot-nlps}
\end{figure}
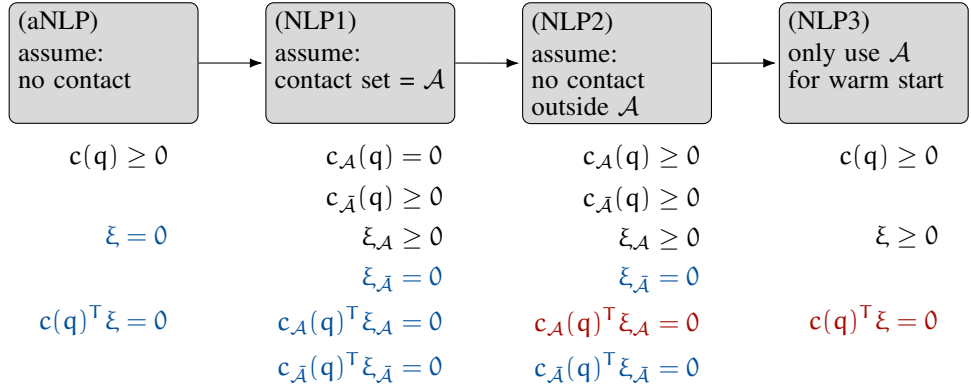

Note that the active set $\setA$ determined in step (0) of our approach
is used both in steps (1) and (2).
In particular, the variables $\xiN$ are implicitly fixed at zero
by omitting them in the implementation.
However, additional contact may occur in steps (1) and (2).
If no solution of \eqref{eq:MPCC} is obtained by solving
\eqref{eq:NLP1} or (NLP2) but the active set differs from $\setA$,
this information is only used in warm starting the variables
but not in choosing the NLP constraints of the following step.
See \cref{fig:quick-shot-nlps} for an overwiew of the four NLPs.

Note also that every solution $x^k$ above must have zero objective value;
otherwise it is a non-physical NLP solution
because the constitutive equation is violated.
We finally note that in the specific examples below
we never obtain solutions of \eqref{eq:MPCC} with (NLP2) or (NLP3):
either the quick shot terminates successfully in step (0) or (1),
or it remains unsuccessful.

\section{Numerical Examples}
\label{sec:example}
In this section, we present numerical examples of the contact problem
with three different obstacles
to show the potential of the proposed quick shot approach.
Specifically, we consider its specialization
to the geometrically exact beam model.
Numerical computations are performed with the elastic
active-set based SQP method developed in \cite{Rose:2018}.
Wherever possible we compute highly accurate NLP solutions
with an SQP tolerance of \num{1e-7};
for the QP subproblem solver (active set method)
we use the default tolerance of $\sqrt{\epsM} \approx \num{1.49e-8}$
where $\epsM \approx \num{2.22e-16}$ is the machine precision.

\begin{figure}
  \centering
  \includegraphics[width=0.5\textwidth]{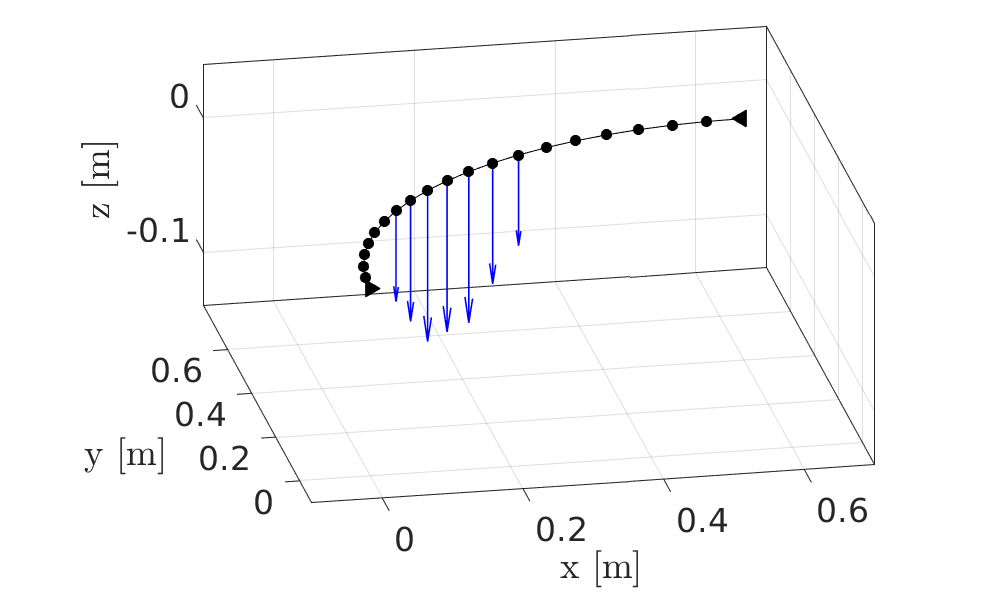}
  \caption{Finite element representation of the beam structure.
    Triangles indicate fixed nodes;
    circles indicate nodes with internal constraints only.
    Arrows denote nodal forces.}
  \label{fig:beam-load}
\end{figure}

All the examples are based on the curved beam structure
that was presented in \cite{Gebhardt:2020a,Gebhardt:2020b}
to illustrate the advantages of the aNLP-based hybrid DDCM framework
in the static and the dynamic case.
Its geometry is described by a quarter of a circular arc
with a total arc length of \qty{1}{m};
see \cref{fig:beam-load}.
Both ends are fully fixed and the structure is uniformly discretized
into 20 two-node finite elements
with a total of 21 nodes numbered $0, \dots, 20$.
No further kinematic restrictions than the internal ones
(orthonormality condition among the three directors) are enforced.
The first node is located at the origin.
The inner nodes can be loaded with arbitrary nodal forces.
\Cref{fig:beam-load} shows the finite element representation
with the external force profile that we apply in this section.
It acts on the inner nodes $7$ to $13$ in the vertical direction
with forces of magnitudes \qty{2}{N}, \qty{4}{N} and \qty{6}{N}.
We will also consider scaled versions of the force profile
with a positive amplification factor $\gamma$.
Of course, our approach can handle any given force profile.

The configuration variable of the resulting NLP is
$q = (q^0, \dots , q^{20}) \in \R^{252}$
with nodal components
$q^i = (\varphi^i_0, d^i_1, d^i_2, d^i_3) \in \R^{12}$.
The strain and stress vectors are $e, s \in \R^{120}$
with six components per finite element each,
similarly the duplicated vectors $\xe, \xs \in \R^{120}$.
Next, there are $24$ boundary conditions fixing $q^0$ and $q^{20}$,
$19 \x 6 = 114$ further kinematic constraints for the directors $d^i_k$
of the inner nodes $i = 1, \dots, 19$,
$120$ compatibility constraints and constitutive equations each,
and finally $19 \x (12 - 6) = 114$ equilibrium constraints
at the inner nodes.
In summary, the relevant dimensions are
$m = 138$, $n = 114$, $n_q = 252$, $n_e = n_s = 120$,
and the NLP has a total of $732$ optimization variables 
with $492$ equality constraints, 
plus a certain number of inequality constraints
that depends on the chosen obstacle.

\newcommand\xshape[3][]{%
  \begin{tikzpicture}[#2]
    \begin{axis}[width=0.6\textwidth,line width=0.3pt,cycle list/Dark2,
      legend pos=north east,legend style={nodes={scale=0.6,transform shape}},
      xmin=0.1,xmax=21.9,ymin=-0.1,ymax=0.02,scaled ticks=false,
      xtick={1,2,3,4,5,6,7,8,9,10,11,12,13,14,15,16,17,18,19,20,21},
      xticklabels={0,,,,,5,,,,,10,,,,,15,,,,,20},
      ytick={-0.2,-0.15,-0.1,-0.05,0,0.05,0.1,0.15},
      yticklabels={$-0.2$,$-0.15$,$-0.1$,$-0.05$,$0$,$0.05$,$0.1$,$0.15$},
      grid,grid style={line width=.1pt,draw=gray!20},
      axis x line=bottom,axis y line=left,#1]
      #3
      \addplot[LUH-blue,thick,mark=square*,mark size=1.5pt]
      coordinates {(1,0) (21,0)};
      \addplot[LUH-blue,only marks,mark=x,semithick] coordinates
      {(2,0) (3,0) (4,0) (5,0) (6,0) (7,0) (8,0) (9,0) (10,0) (11,0)
        (12,0) (13,0) (14,0) (15,0) (16,0) (17,0) (18,0) (19,0) (20,0)};
    \end{axis}
  \end{tikzpicture}
}%
\newcommand\fshape[5][]{\xshape[#1]{}{
    \beamforce{8}{0.02}{$#2$}{red}
    \beamforce{9}{0.04}{$#3$}{red}
    \beamforce{10}{0.06}{}{red}
    \beamforce{11}{0.06}{$#4$}{red}
    \beamforce{12}{0.06}{}{red}
    \beamforce{13}{0.04}{$#3$}{red}
    \beamforce{14}{0.02}{$#2$}{red}
    #5}
}
\newcommand\shape[1][]{\fshape[#1]{}{}{}}

The three obstacles chosen for our contact problem
are illustrated in \cref{tab:obstacles},
where the 3d beam geometry is ``unrolled'' to obtain a simplified 2d view.
\begin{table}
  \centering
  \caption{Obstacles for numerical examples with potential contact nodes.}
  \label{tab:obstacles}
  \def\shape{\xshape[height=3cm,ymajorticks=false]
    {baseline=(current bounding box.center)}}
  \begin{tabular}{lr@{\qquad}c}
    \toprule
    obstacle & nodes & shape (2d view along curved beam) \\
    \midrule
    tip & $10$ &
    \shape{\tipcontact{11}{-0.04}{-0.1}{black!30}} \\
    plane & $4$--$16$ &
    \shape{\levelcontact{5}{17}{-0.04}{-0.1}{black!30}} \\
    hemisphere & $7$--$13$ &
    \shape{\leveldcontact{8}{-0.04}{-0.06}{-0.075}{-0.08}{-0.1}{black!30}} \\
    \bottomrule
  \end{tabular}
\end{table}
The applied load pushes the beam toward the obstacle from above.
To avoid unnecessary technical issues, each obstacle is simply modeled
as a set of lower bounds
on the vertical coordinates of potential contact nodes
rather than a true geometric obstacle.
In terms of the configuration variable $q \in \R^{252}$,
this means that at contact node $i$ we enforce a lower bound on $q_{12 i + 3}$,
the third component of the position vector $\varphi^i_0 \in \R^3$.
The range of node indices $i$ for each obstacle
is given in \cref{tab:obstacles}.

For each of the three obstacles we finally consider
two different material laws,
a symmetric one and an asymmetric one,
both taken from \cite{Gebhardt:2020a,Gebhardt:2020b}.
In terms of the constitutive equation $g(\xe, \xs) = 0$
which defines the constitutive manifold,
the symmetric material law provides
an explicit definition of the stress $\xs$
as a function of the strain $\xe$:
\begin{equation*}
  g(\xe, \xs) = \xs - A \xe - \frac13 A \xe^3.
\end{equation*}
In contrast, the asymmetric material law provides
an implicit definition of the stress $\xs$:
\begin{equation*}
  g(\xe, \xs ) = \xe - \Inv A \xs - 0.0075 \Inv A \xs^2.
\end{equation*}
The matrix $A \in \R^{6 \x 6}$ is diagonal,
$A = \Diag(\SI{75}{N}, \SI{75}{N}, \SI{100}{N},
\SI{100}{Nm^2}, \SI{100}{Nm^2}, \SI{200}{Nm^2})$,
and the powers $\xe^3$, $\xs^2$ are to be understood component-wise.

\subsection{Contact problem 1: tip at node 10}

Here we enforce a single inequality constraint:
the lower bound $q_{123} \ge - 0.1$, 
\cf \cref{tab:obstacles}.
Note that this implies $\setA = \setI = \set{123}$
if the bound becomes active,
and hence (NLP2) and (NLP3) of the quick shot approach become identical
with respective constraints $\cA(q) \ge 0$, $\xiA \ge 0$
and $\cI(q) \ge 0$, $\xi \ge 0$.

\subsubsection{Symmetric manifold}

With the symmetric material law, step (0) of our quick shot approach
does not produce a valid solution of \eqref{eq:MPCC}
for any of the considered values $\gamma \in \set{2,4,6,8,10}$:
solving \eqref{eq:aNLP} always gives an active contact.

\begin{figure}
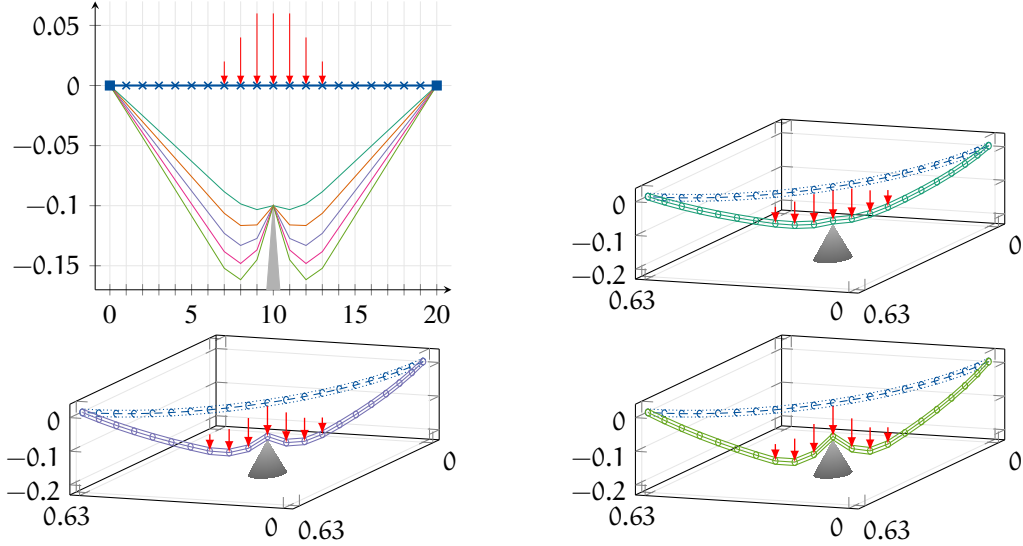

  \centering
  \shape[width=0.42\textwidth,ymin=-0.17,ymax=0.07]
  {\tipcontact{11}{-0.1}{-0.17}{black!30}
    \addplot table[x=node,y index=3]
    {./Dateien/data_exsym_v5_t2_s11_01_opt_p2.dat};
    \addplot table[x=node,y index=3]
    {./Dateien/data_exsym_v5_t4_s11_01_opt_p2.dat};
    \addplot table[x=node,y index=3]
    {./Dateien/data_exsym_v5_t6_s11_01_opt_p2.dat};
    \addplot table[x=node,y index=3]
    {./Dateien/data_exsym_v5_t8_s11_01_opt_p2.dat};
    \addplot table[x=node,y index=3]
    {./Dateien/data_exsym_v5_t10_s11_01_opt_p2.dat};}
  \hfil
  \spike{2}{Dark21}{red}{exsym}{2}
  \\
  \spike{6}{Dark23}{red}{exsym}{2}
  \hfil
  \spike{10}{Dark25}{red}{exsym}{2}
  \caption{Contact problem 1 (tip), symmetric manifold:
    beam deformations computed with step (1) in 2d (top left)
    and in 3d for $\gamma \in \set{2, 6, 10}$
    (top right, bottom left and right).}
  \label{fig:tip-sym}
\end{figure}

In contrast, solving \eqref{eq:NLP1} at step (1),
warm started from the solution of step (0),
immediately yields a valid solution of \eqref{eq:MPCC}
for each value of~$\gamma$.
Thus there is no need to proceed with step (2), equivalently step (3),
but running it for confirmation reproduces the solutions of step (1),
as expected.
In fact, each of the (warm started) confirmation runs
requires just one SQP iteration.
Details of the solutions are given in \cref{tab:tip-sym}
where the measure of the constraint violation is $\norm{c(q)}_\oo$
and the complementarity product is $c(q)\tp \xi$.
The resulting deformations of the beam are illustrated
in \cref{fig:tip-sym}.

\begin{table}
  \centering
  \caption{Contact problem 1 (tip), symmetric manifold:
    solution details of step (1).}
  \label{tab:tip-sym}
  \begin{tabular}{rS[table-format=1.2e+2]
    S[table-format=3.1,table-auto-round]*3{S[table-format=1.2e+1]}cc}
    \toprule
    $\gamma$ & {NLP} & {$\xiA$} & {SQP} & {QP} & {constraint} & compl. & SQP \\
    & {objective}
    & & {tolerance} & {tolerance} & {violation} & product & iter. \\
    \midrule
    $ 2$ & 9.68e-28 & 22.02275591 & \E07 & \tQP & 1.46e-08 & 0 & $3$ \\
    $ 4$ & 4.45e-29 & 73.86633053 & \E07 & \tQP & 1.45e-08 & 0 & $4$ \\
    $ 6$ & 4.04e-32 & 126.0781778 & \E07 & \tQP & 1.33e-08 & 0 & $5$ \\
    $ 8$ & 4.77e-26 & 178.7762652 & \E07 & \tQP & 1.47e-08 & 0 & $5$ \\
    $10$ & 2.68e-23 & 231.9302873 & \E07 & \tQP & 4.59e-08 & 0 & $5$ \\
    \bottomrule
  \end{tabular}
\end{table}

\subsubsection{Asymmetric manifold}

\begin{figure}
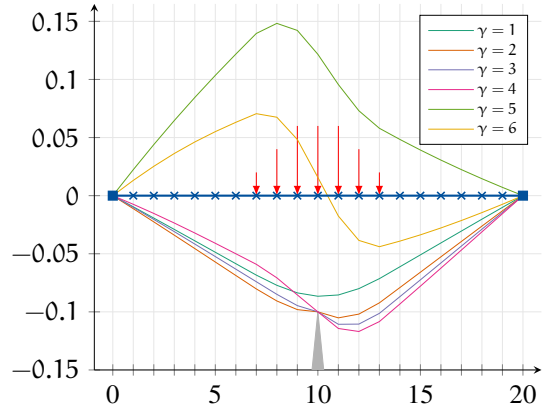

  \centering
  \shape[width=0.5\textwidth,ymin=-0.15,ymax=0.165]
  {\tipcontact{11}{-0.1}{-0.17}{black!30}
    \legend{$\gamma=1$,$\gamma=2$,$\gamma=3$,$\gamma=4$,$\gamma=5$,$\gamma=6$}
    \addplot table[x=node,y index=3] {./Dateien/data_imnon_v5_t1_s11_01.dat};
    \addplot table[x=node,y index=3] {./Dateien/data_imnon_v5_t2_s11_01.dat};
    \addplot table[x=node,y index=3] {./Dateien/data_imnon_v5_t3_s11_01.dat};
    \addplot table[x=node,y index=3] {./Dateien/data_imnon_v5_t4_s11_01.dat};
    \addplot table[x=node,y index=3] {./Dateien/data_imnon_v5_t5_s11_01.dat};
    \addplot table[x=node,y index=3] {./Dateien/data_imnon_v5_t6_s11_01.dat};}
  \caption{Contact problem 1 (tip), asymmetric manifold:
    beam deformations computed with step (0).}
  \label{fig:tip-asym-s0}
\end{figure}

With the asymmetric material law, step (0) of our quick shot approach
converges for all tested values of $\gamma$,
but we need to reduce the SQP and QP solution tolerances for large values.
Moreover, we observe three different types of outcomes
where the active set is either
$\setA = \setI = \set{123}$ (contact)
or $\setA = \0$ (no contact);
see \cref{tab:tip-asym-s0}:
\begin{itemize}
\item For the values $\gamma \in \set{1,5}$
  step (0) produces valid solutions of \eqref{eq:MPCC}
  with zero objective and no contact;
  thus we are successfully finished in these cases.
\item For the values $\gamma \in \set{2,3,4}$
  we obtain valid solutions of (aNLP) but not \eqref{eq:MPCC},
  with contact and nonzero objective (larger than the SQP tolerance);
  thus we proceed with step (1) in these cases.
\item Finally, the value $\gamma = 6$ yields
  again a valid solution of (aNLP),
  but now with nonzero objective and no contact.
  Thus we have a non-physical solution,
  and there are two indications that a physical solution may indeed not exist
  (although we cannot be certain).
  First, the NLP becomes more difficult to solve
  for increasing values of $\gamma$
  and we need to decrease the SQP and QP solution tolerances.
  Second, the beam deformations in \cref{fig:tip-asym-s0}
  clearly show increasingly non-intuitive behavior,
  which becomes quite drastic for the values $\gamma \in \set{5,6}$.
\end{itemize}

\begin{table}
  \centering
  \caption{Contact problem 1 (tip), asymmetric manifold:
    solution details of step (0).}
  \label{tab:tip-asym-s0}
  \def\act{$\set{123}$}
  \sisetup{table-format=1.2e-2,parse-numbers=false}
  \begin{tabular}{cS@{\isep}cSSSc}
    \toprule
    $\gamma$ & {NLP} & {active} & {SQP} & {QP} & {constraint} & SQP \\
    & {objective} & {set $\setA$} & {tolerance}
    & {tolerance} & {violation} & iter. \\
    \midrule
    $1$ & 6.33e-28 & $\0$ & \E07 & \tQP & 2.95e-08 & $3$ \\
    $2$ & 1.28e-01 & \act & \E07 & \tQP & 1.45e-08 & $4$ \\
    $3$ & 3.11e-01 & \act & \E07 & \tQP & 1.46e-08 & $6$ \\
    $4$ & 1.68e-01 & \act & \E06 & \E07 & 9.40e-08 & $5$ \\
    $5$ & 3.06e-17 & $\0$ & \E06 & \E07 & 9.99e-08 & $8$ \\
    $6$ & 2.72e-01 & $\0$ & \E05 & \E06 & 9.96e-07 & $7$ \\
    \bottomrule
  \end{tabular}
\end{table}

\begin{figure}
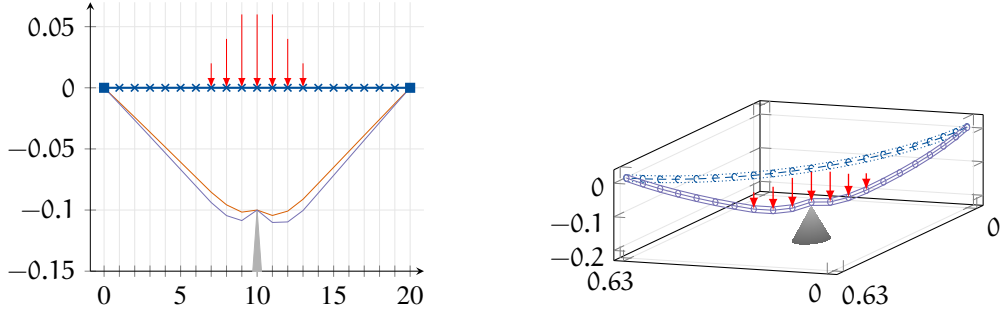

  \centering
  \shape[width=0.4\textwidth,ymin=-0.15,ymax=0.07]
  {\tipcontact{11}{-0.1}{-0.17}{black!30}
    \addplot[Dark22] table[x=node,y index=3]
    {./Dateien/data_imnon_v5_t2_s11_01_opt_p2.dat};
    \addplot[Dark23] table[x=node,y index=3]
    {./Dateien/data_imnon_v5_t3_s11_01_opt_p2.dat};}
  \hfil
  \spike{3}{Dark23}{red}{imnon}{2}
  \caption{Contact problem 1 (tip), asymmetric manifold:
    beam deformations computed with step (1) in 2d (left)
    and for $\gamma = 3$ in 3d (right).}
  \label{fig:tip-asym-s1}
\end{figure}

Solving \eqref{eq:NLP1} in step (1) produces valid solutions of \eqref{eq:MPCC}
for the values $\gamma = 2$ and $\gamma = 3$
whereas the value $\gamma = 4$ yields a positive objective value;
see \cref{tab:tip-asym}.
In comparison to the approximate solutions in \cref{fig:tip-asym-s0},
the asymmetry of the beam deformations for $\gamma \in \set{2,3}$
is less pronounced for the true solutions in \cref{fig:tip-asym-s1}.
As in the symmetric case, running step (2) for confirmation
reproduces these true solutions.
For $\gamma = 4$, on the other hand,
step (2) does not give a valid solution:
the objective value remains positive,
and thus the entire quick shot approach does not succeed.
Details of the optimization runs for this case are listed
in \cref{tab:tip-asym}.
The solutions of both steps are identical up to roundoff errors.

\begin{table}
  \centering
  \caption{Contact problem 1 (tip), asymmetric manifold:
    solution details of step (1) and (if unsuccessful) of step (2).}
  \label{tab:tip-asym}
  \sisetup{table-format=1.2e-2,parse-numbers=false}
  \begin{tabular}{ccS@{\isep}
      S[table-format=2.1,parse-numbers=true,table-auto-round]SSScc}
    \toprule
    $\gamma$
    & step & {NLP} & {$\xiA$} & {SQP} & {QP} & {constraint} & compl. & SQP \\
    & & {objective}
    & & {tolerance} & {tolerance} & {violation} & product & iter. \\
    \midrule
    $2$ & (1) & 1.88e-30 & 21.64426859 & \E07 & \tQP & 1.49e-08 & 0 & $4$ \\
    $3$ & (1) & 2.71e-27 & 47.41802374 & \E07 & \tQP & 1.48e-08 & 0 & $7$ \\
    $4$ & (1) & 1.68e-01 & 0 & \E07 & \tQP & 1.41e-08 & 0 & $2$ \\
	& (2) & 1.68e-01 & 0 & \E07 & \tQP & 1.45e-08 & 0 & $2$ \\
    \bottomrule
  \end{tabular}
\end{table}

\subsection{Contact problem 2: plane at nodes 4--16}

In this example up to 13 inequality constraints can become active:
one at each node from 4 to 16; see also \cref{tab:obstacles}.
Indeed, with the symmetric material law
at least one lower bound becomes active at step (0)
for every chosen value of the factor $\gamma$.
The same observation holds for the asymmetric material law,
except that no contact occurs in the case $\gamma = 6$.
This is apparent from the respective third columns of
\cref{tab:plane-sym-s0,tab:plane-asym-s0},
where the active set is now simply given as the set of contact nodes
(rather than components of the configuration variable~$q$).
Since step (0) does not yield a valid solution of \eqref{eq:MPCC}
in any of the considered cases, we proceed with
steps (1) to (3) of the quick shot approach in case of contact.
\Cref{tab:plane-steps} shows
whether these subsequent steps produce solutions of \eqref{eq:MPCC}.
Further details are discussed in the following.

\begin{table}
  \centering
  \caption{Contact problem 2 (plane), symmetric manifold:
    solution details of step (0).}
  \label{tab:plane-sym-s0}
  \sisetup{table-format=1.2e-2,parse-numbers=false}
  \begin{tabular}{rS@{\isep}cSSSc}
    \toprule
    $\gamma$ & {NLP} & active & {SQP} & {QP} & {constraint} & SQP \\
    & {objective} & set $\setA$
    & {tolerance} & {tolerance} & {violation} & iter. \\
    \midrule
    $ 1$ & 5.55e-02 & $\range{9}{11}$ & \E07 & \tQP & 1.00e-08 & $3$ \\
    $ 2$ & 5.82e-01 & $\range{8}{12}$ & \E07 & \tQP & 1.47e-08 & $4$ \\
    $ 4$ & 2.77e+00 & $\range{8}{12}$ & \E06 & \E07 & 1.00e-07 & $4$ \\
    $ 6$ & 5.67e+00 & $\range{8}{12}$ & \E06 & \E07 & 9.68e-08 & $5$ \\
    $ 8$ & 8.81e+00 & $\range{8}{12}$ & \E07 & \tQP & 1.47e-08 & $5$ \\
    $10$ & 1.20e+01 & $\range{8}{12}$ & \E05 & \E06 & 9.36e-07 & $5$ \\
    \bottomrule
  \end{tabular}
\end{table}

\begin{table}
  \centering
  \caption{Contact problem 2 (plane): overall results of steps (1) to (3).}
  \label{tab:plane-steps}
  \begin{tabular}{rccc}
    \toprule
    $\gamma$ & NLP1 & NLP2 & NLP3 \\
    \midrule
    \mc4c{symmetric manifold} \\
    \midrule
    $1, 2,  4$ & \ok & \ok & \ok \\
    $6, 8, 10$ & \badobj & \badcc & \badcc \\
    \midrule
    \mc4c{asymmetric manifold} \\
    \midrule
    $1$ & \ok & \ok & \ok \\
    $2$ & \badobj & \badobj, \badcc & \badcc \\
    $3$ & \badobj & \badcc & \badcc \\
    $4, 5$ & \badobj & \badobj & \badobj \\
    \bottomrule
  \end{tabular}
\end{table}

\subsubsection{Symmetric manifold}

\Cref{tab:plane-sym-s0} lists solution details of step (0)
for all selected values of the amplification factor $\gamma$.
Overall, there is the tendency that (aNLP) becomes more
difficult to solve with increasing values of $\gamma$,
and in several cases we are only able to obtain solutions
with reduced SQP and QP tolerances.
However, no adjustments are needed for $\gamma = 8$
although it is a relatively large value.

\begin{figure}
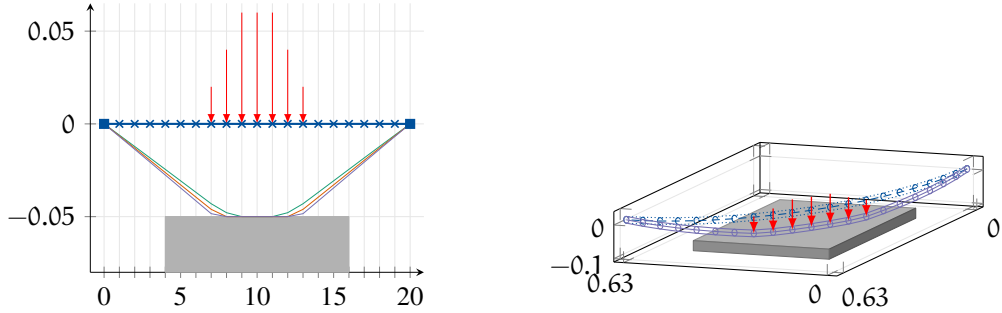

  \centering
  \shape[width=0.4\textwidth,ymin=-0.08,ymax=0.065]
  {\levelcontact{5}{17}{-0.05}{-0.08}{black!30}
    \addplot table[x=node, y=z]
    {./Dateien/data_exsym_v5_t1_e0517_005_opt_p2.dat};
    \addplot table[x=node, y=z]
    {./Dateien/data_exsym_v5_t2_e0517_005_opt_p2.dat};
    \addplot table[x=node, y=z]
    {./Dateien/data_exsym_v5_t4_e0517_005_opt_p2.dat};}
  \hfil
  \plane{4}{Dark23}{red}{exsym}{2}
  \caption{Contact problem 2 (plane), symmetric manifold:
    beam deformation computed with step (1) in 2d (left)
    and for $\gamma = 4$ in 3d (right).}
  \label{fig:plane-sym}
\end{figure}

With force amplifications $\gamma \in \set{1, 2, 4}$,
the respective solutions of \eqref{eq:NLP1} have zero objective,
and hence we obtain valid solutions of \eqref{eq:MPCC}.
\Cref{tab:plane-sym-s1} gives the solution details for these cases,
and \cref{fig:plane-sym} illustrates the corresponding beam deformations.
As in contact problem~1, running the unnecessary steps (2) and (3)
for confirmation reproduces these valid solutions.

\begin{table}
  \centering
  \caption{Contact problem 2 (plane), symmetric manifold:
    solution details of step (1) in successful cases.}
  \label{tab:plane-sym-s1}
  \sisetup{table-format=1.2e-2,parse-numbers=false}
  \begin{tabular}{cS@{\isep}cSScc}
    \toprule
    $\gamma$ & {NLP} & {$\xiA$} & {SQP} & {constraint} & compl. & SQP \\
    & {objective} & & {tolerance} & {violation} & product & iter. \\
    \midrule
    $1$ & 3.30e-29
    & $(       2.86, 6.06, 2.86      )$ & \E07 & 1.48e-08 & 0 & $3$ \\
    $2$ & 6.98e-30
    & $( 2.15,12.06,12.06,12.06, 2.15)$ & \E07 & 1.49e-08 & 0 & $3$ \\
    $4$ & 5.57e-36
    & $(13.65,24.06,24.06,24.06,13.65)$ & \E07 & 1.42e-08 & 0 & $4$ \\
    \bottomrule
  \end{tabular}
\end{table}

For the remaining values $\gamma \in \set{6, 8, 10}$,
\eqref{eq:NLP1} yields non-physical solutions with positive objective values.
In the subsequent steps (2) and (3), each objective value
drops to zero within the respective SQP tolerance,
which needs to be slightly reduced in step (2) for $\gamma = 8$
and in step (1) for $\gamma = 10$.
However, the complementarity products become nonzero
so that we still do not arrive at valid solutions of \eqref{eq:MPCC}.
To verify these results, we repeat step (3)
for $\gamma \in \set{6, 8, 10}$ with higher accuracy,
obtaining solutions with SQP tolerances up to \num{1e-9} or \num{1e-10}
with respective QP tolerances of \num{1e-10} or \num{1e-11}.
Again all objectives are zero within the given tolerances
whereas the complementarity products remain nonzero.
Detailed solution statistics of these unsuccessful cases are provided
in \cref{tab:plane-sym}
for steps (1) to (3) of the quick shot approach,
including the verification runs with increased accuracy.

\begin{table}
  \centering
  \caption{Contact problem 2 (plane), symmetric manifold:
    solution details of steps (1) to (3)
    in cases where step (1) is unsuccessful.}
  \label{tab:plane-sym}
  \sisetup{table-format=1.2e-2,parse-numbers=false}
  \begin{tabular}{rcSSSSS[table-format=1.3]S[table-format=3]}
    \toprule
    $\gamma$ & step & {NLP} & {SQP} & {QP} & {constraint} & {compl.} & {SQP} \\
    & & {objective}
    & {tolerance} & {tolerance} & {violation} & {product} & {iter.} \\
    \midrule
    $ 6$ & (1) & 2.40e-04 & \E07 & \tQP & 1.45e-08 & 0     &   4 \\
         & (2) & 2.44e-09 & \E07 & \tQP & 1.48e-08 & 0.135 &  69 \\
         & (3) & 1.54e-09 & \E07 & \tQP & 1.47e-08 & 0.135 &   3 \\
         & (3) & 4.63e-15 & \E10 & \E11 & 9.81e-12 & 0.135 & 108 \\
    \addlinespace
    $ 8$ & (1) & 4.25e-03 & \E07 & \tQP & 1.40e-08 & 0     &   4 \\
         & (2) & 2.23e-07 & \E06 & \E07 & 9.72e-08 & 0.781 &  57 \\
         & (3) & 1.60e-09 & \E07 & \tQP & 1.48e-08 & 0.782 &  32 \\
         & (3) & 4.55e-15 & \E10 & \E11 & 1.00e-11 & 0.782 & 155 \\
    \addlinespace
    $10$ & (1) & 1.31e-02 & \E06 & \E07 & 9.70e-08 & 0     &   4 \\
         & (2) & 2.31e-09 & \E07 & \tQP & 1.48e-08 & 1.752 &  99 \\
         & (3) & 1.83e-09 & \E07 & \tQP & 1.48e-08 & 1.127 &  97 \\
         & (3) & 1.12e-12 & \E09 & \E10 & 9.98e-11 & 1.127 & 239 \\
    \bottomrule
  \end{tabular}
\end{table}

In this example we observe that the active sets
computed in steps (1) to (3) differ
from the active set of step (0) for the values $\gamma \in \set{6, 8, 10}$.
\cref{tab:plane-sym-A} lists the active sets encountered in each step.
We finally notice that similar behavior also occurs
with even larger amplifications of the force profile.

\begin{table}
  \centering
  \caption{Contact problem 2 (plane), symmetric manifold: active sets $\setA$.}
  \label{tab:plane-sym-A}
  \begin{tabular}{ccc@{\qquad}ccc@{\qquad}ccc}
    \toprule
    $\gamma$ & step & $\setA$ &
    $\gamma$ & step & $\setA$ &
    $\gamma$ & step & $\setA$ \\
    \midrule
    $ 6$ & (0) & $\range{8}{12}$ &
    $ 8$ & (0) & $\range{8}{12}$ &
    $10$ & (0) & $\range{8}{12}$ \\
         & (1) & $\range{7}{13}$ &
         & (1) & $\range{7}{13}$ &
         & (1) & $\range{7}{13}$ \\
         & (2) & $\set{7,13}$ &
         & (2) & $\set{7,13}$ &
         & (2) & $\set{7,13}$ \\
         & (3) & $\set{7,13}$ &
         & (3) & $\set{7,13}$ &
         & (3) & $\set{10}$ \\
    \bottomrule
  \end{tabular}
\end{table}

\subsubsection{Asymmetric manifold}

\begin{figure}
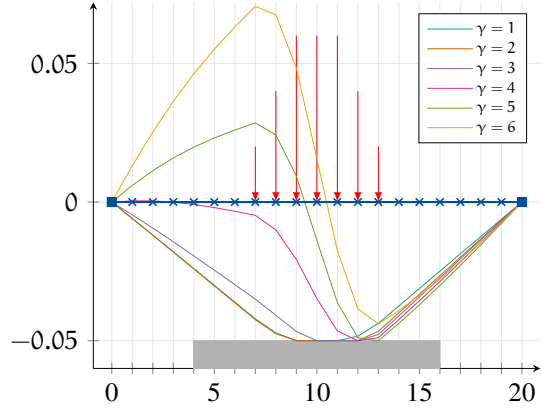

  \centering
  \shape[width=0.5\textwidth,ymin=-0.06,ymax=0.072]
  {\levelcontact{5}{17}{-0.05}{-0.06}{black!30}
    \legend{$\gamma=1$,$\gamma=2$,$\gamma=3$,$\gamma=4$,$\gamma=5$,$\gamma=6$}
    \addplot table[x=node,y=z] {./Dateien/data_imnon_v5_t1_e0517_005.dat};
    \addplot table[x=node,y=z] {./Dateien/data_imnon_v5_t2_e0517_005.dat};
    \addplot table[x=node,y=z] {./Dateien/data_imnon_v5_t3_e0517_005.dat};
    \addplot table[x=node,y=z] {./Dateien/data_imnon_v5_t4_e0517_005.dat};
    \addplot table[x=node,y=z] {./Dateien/data_imnon_v5_t5_e0517_005.dat};
    \addplot table[x=node,y=z] {./Dateien/data_imnon_v5_t6_e0517_005.dat};}
  \caption{Contact problem 2 (plane), asymmetric manifold:
    beam deformations computed with step (0).}
  \label{fig:plane-asym-s0}
\end{figure}

With the asymmetric, implicit material law,
the solutions of step (0) show again non-intuitive behavior:
as expected, the contact set becomes initially larger with increasing force,
but already at $\gamma = 3$ the growth reverts to shrinking.
These results are similar to those when the tip is chosen as an obstacle
(\cref{fig:tip-asym-s0}).
All the NLP solutions obtained at step (0) are non-physical
with positive objective value, see \cref{tab:plane-asym-s0}.
The computed (non-physical) beam deformations
are illustrated in \cref{fig:plane-asym-s0}.
Once again the problem becomes more difficult to solve
with increasing values of $\gamma$,
and we need to decrease the SQP and QP solution tolerances
for $\gamma \in \set{4, 5, 6}$.
For the values $\gamma = 5$ and $\gamma = 6$,
we are not even able to solve \eqref{eq:aNLP} from scratch.
In fact, we use the solution for $\gamma = 4$
as a warm start to obtain the solution for $\gamma = 5$,
and likewise with $\gamma = 5$ and $\gamma = 6$.
Despite the warm start, the SQP tolerance
has to be reduced to \num{1e-4} in the two last cases.
\begin{table}
  \centering
  \caption{Contact problem 2 (plane), asymmetric manifold: results of step (0).
    Iterations marked with * are warm started
    from the solution for $\gamma - 1$.}
  \label{tab:plane-asym-s0}
  \sisetup{table-format=1.2e-2,parse-numbers=false}
  \def\1{\hphantom1}
  \begin{tabular}{cS@{\isep}cSSSc}
    \toprule
    $\gamma$ & {NLP} & active & {SQP} & {QP} & {constraint} & SQP \\
    & {objective} & set $\setA$
    & {tolerance} & {tolerance} & {violation} & iter. \\
    \midrule
    $1$ & 5.16e-02 & $\range{ 9}{11}$ & \E07 & \tQP & 1.47e-08 & $4$ \\
    $2$ & 4.34e-01 & $\range{ 9}{12}$ & \E07 & \tQP & 1.49e-08 & $6$ \\
    $3$ & 6.50e-01 & $\range{10}{12}$ & \E07 & \tQP & 1.48e-08 & $8$ \\
    $4$ & 3.09e-01 & $\set{12}$       & \E06 & \E07 & 9.98e-08 & $6$ \\
    $5$ & 6.20e-03 & $\set{13}$       & \E04 & \E05 & 1.00e-05 & $3\rlap*$ \\
    $6$ & 2.72e-01 & $\0$             & \E04 & \E05 & 1.62e-05 & $6\rlap*$ \\
    \bottomrule
  \end{tabular}
\end{table}

\begin{figure}
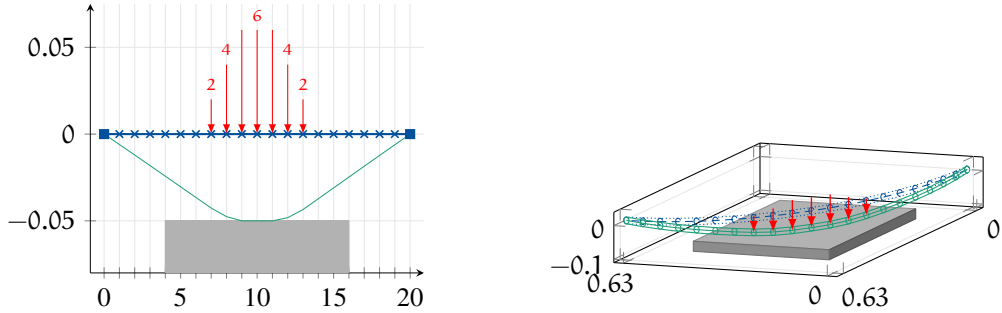

  \centering
  \fshape[width=0.4\textwidth,ymin=-0.08,ymax=0.075]{2}{4}{6}
  {\levelcontact{5}{17}{-0.05}{-0.08}{black!30}
    \addplot table[x=node,y=z]
    {./Dateien/data_imnon_v5_t1_e0517_005_opt_p2.dat};}
  \hfil
  \plane{1}{Dark21}{red}{imnon}{2}
  \caption{Contact problem 2 (plane), asymmetric manifold, $\gamma = 1$:
    beam deformation computed with step (1).}
  \label{fig:plane-asym-s1}
\end{figure}

In these examples, the solution of step (1) of our approach
has zero objective value only for $\gamma = 1$,
so that a valid solution of \eqref{eq:MPCC} is found.
Details of this solution are given in \cref{tab:plane-asym-s1},
\cref{fig:plane-asym-s1} shows the resulting beam deformation.
Again, solving (NLP2) and (NLP3) for confirmation
reproduces the solution of step (1).
\begin{table}
  \centering
  \caption{Contact problem 2 (plane), asymmetric manifold, $\gamma = 1$:
    valid solution of step (1) for \eqref{eq:MPCC}.}
  \label{tab:plane-asym-s1}
  \sisetup{table-format=1.2e-2,parse-numbers=false}
  \begin{tabular}{cS@{\isep}cSSScc}
    \toprule
    $\gamma$ & {NLP} & $\xi$ & {SQP} & {QP} & {constraint} & {compl.} & SQP \\
    & {objective}
    & & {tolerance} & {tolerance} & {violation} & {product} & iter. \\
    \midrule
    $1$ & 2.40e-27 & $(2.29,6.07,3.39)$ & \E07 & \tQP & 1.46e-08 & 0 & $3$ \\
    \bottomrule
  \end{tabular}
\end{table}

All other values of $\gamma$ yield
non-physical numerical solutions in step (1),
with positive objective values.
For the amplification factors $\gamma = 4$ and $\gamma = 5$,
these non-physical solutions are reproduced in steps (2) and (3),
and all of them satisfy the complementarity condition.
For the amplification factors $\gamma = 2$ and $\gamma = 3$, in contrast,
the objective value drops to zero (or almost to zero) in steps (2) and (3)
whereas the complementarity product becomes positive.
Thus our quickshot approach does not succeed
for any value $\gamma \in \set{2, 3, 4, 5}$.
Again we need to reduce the SQP and QP tolerances in several cases,
and as in the symmetric case we verify the results
by repeating step (3) with increased accuracy.
Here we are able to obtain solutions with SQP and QP tolerances
up to \num{1e-10} and \num{1e-11}, respectively.
Detailed solution statistics of steps (1) to (3) are given in
\Cref{tab:plane-asym} for $\gamma \in \set{2, 3, 4, 5}$,
again including the highly accurate verification runs.

\begin{table}
  \centering
  \caption{Contact problem 2 (plane), asymmetric manifold:
    solution details of steps (1) to (3)
    in cases where step (1) is unsuccessful.}
  \label{tab:plane-asym}
  \sisetup{table-format=1.2e-2,parse-numbers=false}
  \begin{tabular}{rcSSSS@{\isep}
    S[table-format=+1.2e-2,parse-numbers=true,
    round-mode=figures,round-precision=3,round-pad=false]
    S[table-format=3]}
    \toprule
    $\gamma$ & step & {NLP} & {SQP} & {QP} & {constraint} & {compl.} & {SQP} \\
    & & {objective}
    & {tolerance} & {tolerance} & {violation} & {product} & {iter.} \\
    \midrule
    $2$ & (1) & 1.59e-04 & \E07 & \tQP & 1.48e-08 & -1.7572e-15 &   5 \\
        & (2) & 1.79e-05 & \E05 & \E06 & 9.85e-07 &  0.0304     &   6 \\
        & (3) & 1.21e-09 & \E07 & \tQP & 1.49e-08 &  0.0326     &  79 \\
        & (3) & 3.88e-13 & \E09 & \E10 & 9.97e-11 &  0.0326     & 156 \\
    \addlinespace
    $3$ & (1) & 6.00e-02 & \E07 & \tQP & 1.48e-08 & 1.2924e-13  & 207 \\
        & (2) & 1.60e-26 & \E05 & \E06 & 9.85e-07 & 1.7674      &   3 \\
        & (3) & 5.85e-12 & \E07 & \tQP & 1.41e-08 & 1.7674      &   1 \\
        & (3) & 1.13e-14 & \E10 & \E11 & 8.15e-12 & 1.7674      & 140 \\
    \addlinespace
    $4$ & (1) & 3.09e-01 & \E07 & \tQP & 1.46e-08 & 0           &   2 \\
        & (2) & 3.09e-01 & \E07 & \tQP & 1.44e-08 & 0           &   2 \\
        & (3) & 3.09e-01 & \E07 & \tQP & 1.43e-08 & 0           &   2 \\
        & (3) & 3.09e-01 & \E10 & \E11 & 9.50e-12 & 0           &   3 \\
    \addlinespace
    $5$ & (1) & 6.19e-03 & \E07 & \tQP & 1.36e-08 & 0           &   2 \\
        & (2) & 6.19e-03 & \E06 & \E07 & 9.43e-08 & 0           &   1 \\
        & (3) & 6.19e-03 & \E06 & \E07 & 9.99e-08 & 0           &   1 \\
        & (3) & 6.19e-03 & \E10 & \E11 & 9.43e-12 & 0           &   3 \\
    \bottomrule
  \end{tabular}
\end{table}

We also observe that the active sets
computed in steps (1) to (3) differ
from the active set of step (0) for the values $\gamma \in \set{2,3,4,5}$.
\cref{tab:plane-asym-A} lists the active sets encountered in each step.
We finally notice that similar behavior also occurs
with even larger amplifications of the force profile.

\begin{table}
  \centering
  \caption{Contact problem 2 (plane), asymmetric manifold: active sets $\setA$.}
  \label{tab:plane-asym-A}
  \begin{tabular}{ccc@{\qquad}ccc@{\qquad}ccc@{\qquad}ccc}
    \toprule
    $\gamma$ & step & $\setA$ &
    $\gamma$ & step & $\setA$ &
    $\gamma$ & step & $\setA$ &
    $\gamma$ & step & $\setA$ \\
    \midrule
    $ 2$ & (0) & $\range{9}{12}$ &
    $ 3$ & (0) & $\range{10}{12}$ &
    $ 4$ & (0) & $\set{12}$ &
    $ 5$ & (0) & $\set{13}$ \\
         & (1) & $\range{8}{12}$ &
         & (1) & $\range{8}{12}$ &
         & (1) & $\set{12}$ &
         & (1) & $\set{13}$ \\
         & (2) & $\set{8}$ &
         & (2) & $\set{8}$ &
         & (2) & $\set{12}$ &
         & (2) & $\set{13}$ \\
         & (3) & $\set{8}$ &
         & (3) & $\set{8}$ &
         & (3) & $\set{12}$ &
         & (3) & $\set{13}$ \\
    \bottomrule
  \end{tabular}
\end{table}

\subsection{Contact problem 3: convex hemisphere at nodes 7--13}

In this example up to seven inequality constraints can become active:
the lower bounds at nodes $7$ to $13$, see \cref{tab:obstacles}.
\Cref{tab:hemis-steps} gives an overview of all results
obtained with our quick shot approch on this problem.
Here we omit the cases $\gamma = 1$ and $\gamma = 2$
where step (0) already produces valid solutions of \eqref{eq:MPCC}
with both the symmetric and the asymmetric material laws.
\begin{table}
  \centering
  \caption{Contact problem 3 (hemisphere): overall results of steps (1) to (3).
  NLP2 and NLP3 coincide in row two of the explicit case since $\setA = \setI$.}
  \label{tab:hemis-steps}
  \begin{tabular}{rccc}
    \toprule
    {$\gamma$} & NLP1 & NLP2 & NLP3 \\
    \midrule
    \mc4c{symmetric manifold} \\
    \midrule
    2.25, 2.5 & \ok & \ok & \ok \\
    2.75, 3, 4, 6, 8, 10 & \ok & \mc2c{\ok} \\
    \midrule
    \mc4c{asymmetric manifold} \\
    \midrule
    2.5, 3, 2.5 & \badobj & \badcc & \badcc \\
    4 & \badobj & diverged & \badcc \\
    \bottomrule
  \end{tabular}
\end{table}

\subsubsection{Symmetric manifold}

\begin{figure}
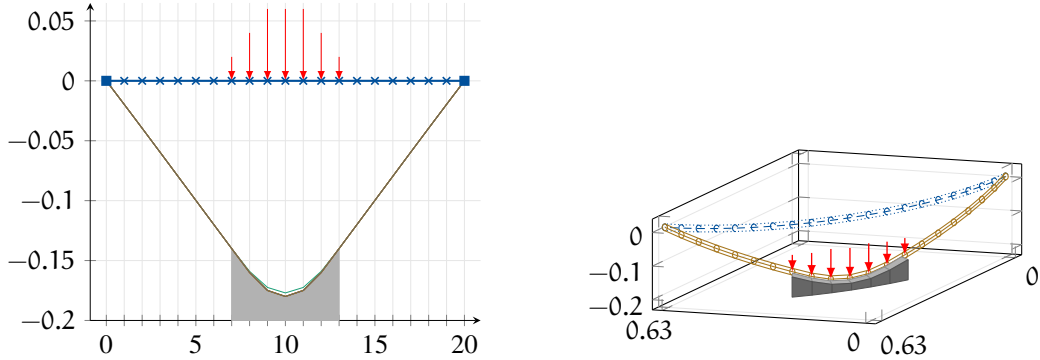

  \centering
  \shape[width=0.45\textwidth,ymin=-0.2,ymax=0.065]
  {\leveldcontact{8}{-0.14}{-0.16}{-0.175}{-0.18}{-0.2}{black!30}
    \addplot table[x=node,y=z]
    {./Dateien/data_exsym_v5_t225_ed0814_018014_opt_p2.dat};
    \addplot table[x=node,y=z]
    {./Dateien/data_exsym_v5_t25_ed0814_018014_opt_p2.dat};
    \addplot table[x=node,y=z]
    {./Dateien/data_exsym_v5_t275_ed0814_018014_opt_p2.dat};
    \addplot table[x=node,y=z]
    {./Dateien/data_exsym_v5_t3_ed0814_018014_opt_p2.dat};
    \addplot table[x=node,y=z]
    {./Dateien/data_exsym_v5_t4_ed0814_018014_opt_p2.dat};
    \addplot table[x=node,y=z]
    {./Dateien/data_exsym_v5_t6_ed0814_018014_opt_p2.dat};
    \addplot table[x=node,y=z]
    {./Dateien/data_exsym_v5_t8_ed0814_018014_opt_p2.dat};
    \addplot table[x=node,y=z]
    {./Dateien/data_exsym_v5_t10_ed0814_018014_opt_p2.dat};}
  \hfil
  \hemis{10}{Dark27}{red}{exsym}{2}
  \caption{Contact problem 3 (hemisphere), symmetric manifold:
    beam deformations computed with step (1) in 2d (left)
    and for $\gamma = 10$ in 3d (right).}
  \label{fig:hemis-sym}
\end{figure}

For the symmetric material law, we report results on force amplification
factors $\gamma$ between $1$ and $10$.
For the two smallest values $\gamma = 1$ and $\gamma = 2$,
no contact occurs and step (0) of our quick shot approach
generates valid solutions of \eqref{eq:MPCC}.
For $\gamma \ge 3$ we obtain non-physical solutions
with positive objective value
and with all seven bounds being active,
which means that $\setA = \setI$ and (NLP2) is identical to (NLP3).
To investigate the transition range, we also consider
the values $\gamma \in \set{2.25, 2.5, 2.75}$.
It turns out that all seven bounds remain active for $\gamma = 2.75$
whereas we obtain proper subsets $\setA \subset \setI$
for $\gamma \in \set{2.25, 2.5}$.
Solution details of step (0)  for all selected values of $\gamma$
are listed in \cref{tab:hemis-sym-s0}.
\begin{table}
  \centering
  \caption{Contact problem 3 (hemisphere), symmetric manifold:
    solution details of step (0).}
  \label{tab:hemis-sym-s0}
  \sisetup{table-format=1.2e-2,parse-numbers=false}
  \begin{tabular}{S[table-format=2.2]S@{\isep}cSSSc}
    \toprule
    {$\gamma$} & {NLP} & {active} & {SQP} & {QP} & {constraint} & {SQP} \\
    & {objective} & {set $\setA$} & {tolerance}
    & {tolerance} & {violation} & {iter.} \\
    \midrule
    1    & 1.69e-28 & $\0$              & \E07 & \tQP & 1.17e-08 & 3 \\
    2    & 1.78e-32 & $\0$              & \E07 & \tQP & 1.49e-08 & 4 \\
    2.25 & 7.42e-03 & $\set{7,13}$      & \E07 & \tQP & 1.43e-08 & 4 \\
    2.5  & 3.94e-02 & $\set{7,8,12,13}$ & \E07 & \tQP & 1.46e-08 & 4 \\
    2.75 & 9.73e-02 & $\range{7}{13}$   & \E07 & \tQP & 1.44e-08 & 4 \\
    3    & 1.80e-01 & $\range{7}{13}$   & \E07 & \tQP & 1.43e-08 & 4 \\
    4    & 7.16e-01 & $\range{7}{13}$   & \E07 & \tQP & 1.31e-08 & 4 \\
    6    & 2.44e+00 & $\range{7}{13}$   & \E07 & \tQP & 1.40e-08 & 5 \\
    8    & 4.60e+00 & $\range{7}{13}$   & \E07 & \tQP & 1.46e-08 & 5 \\
    10   & 6.98e+00 & $\range{7}{13}$   & \E04 & \E05 & 9.98e-06 & 5 \\
    \bottomrule
  \end{tabular}
\end{table}
As mentioned, step (0) is successful for $\gamma = 1$ and $\gamma = 2$.
For all other values of $\gamma$, the numerical solution of \eqref{eq:aNLP}
has positive objective value and is thus not a solution of \eqref{eq:MPCC};
for $\gamma = 10$ we even need to reduce the solution tolerances.
The solution of \eqref{eq:NLP1} in step (1)
has zero objective in all cases with $\gamma > 2$.
Consequently, step (1) produces valid solutions of \eqref{eq:MPCC}
where step (0) does not.
When running steps (2) and (3) for confirmation,
these solutions are reproduced by (NLP2) and (NLP3)
(where in fact steps (2) and (3) are identical
except for $\gamma \in \set{2.25, 2.5}$).
\Cref{fig:hemis-sym} shows the beam deformations
for the solutions with contact.
Interestingly, solutions of (NLP3) in step (3)
are also obtained with significantly higher accuracy:
in all cases we reach an SQP tolerance of \num{1e-10}
with a QP tolerance of \num{1e-11}.
Solution details including the increased tolerances
are given in \cref{tab:hemis-sym}.

\begin{table}
  \centering
  \caption{Contact problem 3 (hemisphere), symmetric manifold:
    solution details of step (1) where successful.}
  \sisetup{table-format=1.2e-2,parse-numbers=false}
  \begin{tabular}{S[table-format=2.2]S@{\isep}cSSc}
    \toprule
    {$\gamma$} & {NLP}& {$\xiA$} & {SQP} & {constraint} & {SQP} \\
    & {objective} & & {tolerance} & {violation} & {iter.} \\
    \midrule
    2.25 & 4.84e-31 &
    $(2.83,2.83)$ & \E07 & 1.46e-08 & $3$ \\
    2.5  & 1.38e-27 &
    $(5.02,1.57,1.57,5.02)$ & \E07 & 1.49e-08 & $3$ \\
    2.75 & 6.43e-28 &
    $(5.52,2.84,1.12,1.71,1.12,2.84,5.62)$ & \E07 & 1.46e-08 & $3$ \\
    3    & 5.88e-30 &
    $(6.02,3.84,2.62,3.21,2.62,3.84,6.02)$ & \E07 & 1.46e-08 & $3$ \\
    4    & 2.37e-28 &
    $(8.02,7.84,8.62,9.21,8.62,7.84,8.02)$ & \E07 & 1.46e-08 & $3$ \\
    6    & 2.51e-29 &
    $(12.02,15.84,20.62,21.21,20.62,15.84,12.02)$ & \E07 & 1.40e-08 & $4$ \\
    8    & 5.38e-29 &
    $(16.02,23.84,32.62,33.21,32.62,23.84,16.02)$ & \E07 & 1.42e-08 & $4$ \\
    10   & 3.49e-29 &
    $(20.02,31.84,44.62,45.21,44.62,31.84,20.02)$ & \E07 & 1.40e-08 & $4$ \\
    \bottomrule
  \end{tabular}
\end{table}

\begin{table}
  \centering
  \caption{Contact problem 3 (hemisphere), symmetric manifold:
    solution details of steps (1) to (3) where successful.
    Step (3) is identical to step (2) for $\gamma \ge 2.75$.}
  \label{tab:hemis-sym}
  \sisetup{table-format=1.2e-2,parse-numbers=false}
  \begin{tabular}{S[table-format=2.2]cSSSSS[table-format=2]}
    \toprule
    {$\gamma$} & step & {NLP} & {SQP} & {QP} & {constraint}  & {SQP} \\
    & & {objective} & {tolerance} & {tolerance} & {violation} & {iter.} \\
    \midrule
    2.25 & (1) & 4.84e-31 & \E07 & \tQP & 1.46e-08 &  3 \\
         & (2) & 5.32e-14 & \E07 & \tQP & 1.38e-08 &  1 \\
         & (3) & 3.90e-33 & \E07 & \tQP & 1.44e-08 &  1 \\
         & (3) & 6.59e-31 & \E10 & \E11 & 1.34e-11 &  1 \\
    \addlinespace
    2.5  & (1) & 1.38e-27 & \E07 & \tQP & 1.49e-08 &  3 \\
         & (2) & 9.74e-14 & \E07 & \tQP & 1.46e-08 &  1 \\
         & (3) & 4.20e-27 & \E07 & \tQP & 1.38e-08 &  1 \\
         & (3) & 2.06e-31 & \E10 & \E11 & 2.70e-11 &  1 \\
    \addlinespace
    2.75 & (1) & 6.43e-28 & \E07 & \tQP & 1.46e-08 &  3 \\
         & (2) & 1.86e-28 & \E07 & \tQP & 1.47e-08 &  1 \\
         & (2) & 3.47e-15 & \E10 & \E11 & 9.48e-12 & 18 \\
    \addlinespace
    3    & (1) & 5.88e-30 & \E07 & \tQP & 1.46e-08 &  3 \\
         & (2) & 1.61e-13 & \E07 & \tQP & 1.47e-08 &  1 \\
         & (2) & 4.93e-15 & \E10 & \E11 & 9.03e-12 & 35 \\
    \addlinespace
    4    & (1) & 2.37e-28 & \E07 & \tQP & 1.46e-08 &  3 \\
         & (2) & 1.96e-13 & \E07 & \tQP & 1.48e-08 &  1 \\
         & (2) & 6.13e-15 & \E10 & \E11 & 9.87e-12 & 38 \\
    \addlinespace
    6    & (1) & 2.51e-29 & \E07 & \tQP & 1.40e-08 &  4 \\
         & (2) & 4.91e-14 & \E07 & \tQP & 1.47e-08 &  1 \\
         & (2) & 5.38e-15 & \E10 & \E11 & 9.48e-12 & 25 \\
    \addlinespace
    8    & (1) & 5.38e-29 & \E07 & \tQP & 1.42e-08 &  4 \\
         & (2) & 5.93e-14 & \E07 & \tQP & 1.47e-08 &  1 \\
         & (2) & 4.81e-15 & \E10 & \E11 & 9.96e-12 & 27 \\
    \addlinespace
    10   & (1) & 3.49e-29 & \E07 & \tQP & 1.40e-08 &  4 \\
         & (2) & 4.43e-14 & \E07 & \tQP & 1.47e-08 &  1 \\
         & (2) & 6.81e-15 & \E10 & \E11 & 9.50e-12 & 17 \\
    \bottomrule
  \end{tabular}
\end{table}

\subsubsection{Asymmetric manifold}

\begin{figure}
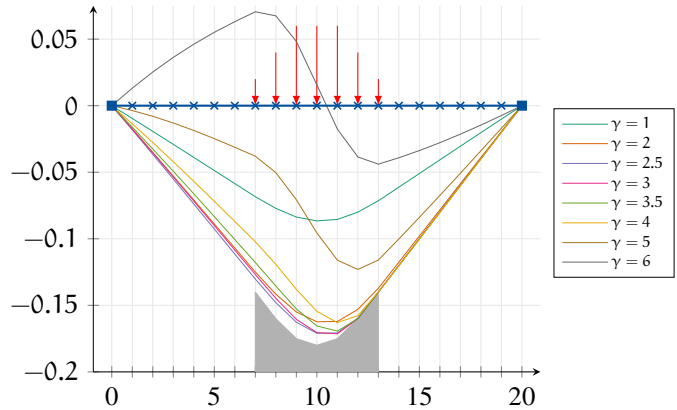

  \centering
  \shape[width=0.5\textwidth,ymin=-0.2,ymax=0.075,
  legend style={cells={anchor=west},at={(1.3,0.72)}}]
  {\leveldcontact{8}{-0.14}{-0.16}{-0.175}{-0.18}{-0.2}{black!30}
    \legend{$\gamma=1$,$\gamma=2$,$\gamma=2.5$,$\gamma=3$,
      $\gamma=3.5$,$\gamma=4$,$\gamma=5$,$\gamma=6$}
    \addplot table[x=node,y index=3]
    {./Dateien/data_imnon_v5_t1_ed0814_018014.dat};
    \addplot table[x=node,y index=3]
    {./Dateien/data_imnon_v5_t2_ed0814_018014.dat};
    \addplot table[x=node,y index=3]
    {./Dateien/data_imnon_v5_t25_ed0814_018014.dat};
    \addplot table[x=node,y index=3]
    {./Dateien/data_imnon_v5_t3_ed0814_018014.dat};
    \addplot table[x=node,y index=3]
    {./Dateien/data_imnon_v5_t35_ed0814_018014.dat};
    \addplot table[x=node,y index=3]
    {./Dateien/data_imnon_v5_t4_ed0814_018014.dat};
    \addplot table[x=node,y index=3]
    {./Dateien/data_imnon_v5_t5_ed0814_018014.dat};
    \addplot table[x=node,y index=3]
    {./Dateien/data_imnon_v5_t6_ed0814_018014.dat};}
  \caption{Contact problem 3 (hemisphere), asymmetric manifold:
    beam deformations computed with step (0).}
  \label{fig:hemis-asym-s0}
\end{figure}

With the asymmetric implicit material law,
the solutions of \eqref{eq:aNLP} at step (0)
show the same non-intuitive behavior as with the tip and the plane:
when increasing the forces, contact occurs at intermediate values,
here $\gamma \in \set{2.5, 3, 3.5, 4}$,
whereas smaller and larger values do not lead to contact,
see \cref{tab:hemis-asym-s0} and \cref{fig:hemis-asym-s0}.
For the values 
where contact occurs in step (0),
further computations are performed.
Unfortunately, a solution of \eqref{eq:MPCC}
cannot be found for any of these cases with our quick shot approach.
In step (1), all numerical solutions have positive objective value,
and for the values $\gamma \in \set{3, 3.5, 4}$
the SQP and QP tolerances need to be reduced to
\num{1e-6} and \num{1e-7}, repectively, see \cref{tab:hemis-asym}.
Next, in step (2), all numerical solutions have zero objective
but violate the complementarity condition,
and again we need to reduce the tolerances.
Furthermore, for $\gamma = 4$,
the SQP algorithm does not even converge in step (2)
so that no numerical solution is available.
In the subsequent step (3), the warm start is then performed
with $(x^1, \xi^1)$ from step (1),
and a solution is computed with substantially reduced tolerances
of \num{1e-3} for the SQP method and \num{1e-4} for the QP solver.
In all other cases, i.e., for $\gamma < 4$,
step (3) produces numerical solutions with zero objective
which still violate the complementarity condition.
Moreover, in all these cases step (3) is successfully repeated
with significantly increased tolerances
of \num{1e-10} (SQP) and \num{1e-11} (QP).
However, this just confirms the violation of the complementarity condition.

Interestingly, we obtain a solution of \eqref{eq:MPCC} with step (0)
not only for the values $\gamma \in \set{1, 2}$ but also for $\gamma = 5$.
Recall that this means a solution without contact and with zero objective.
For larger force amplification values $\gamma \ge 6$,
there is still no contact but the objective values become positive.
In fact, they grow drastically with $\gamma$,
and the numerical difficulties increase as well:
the SQP and QP solution tolerances have to be reduced in most cases,
and for $\gamma = 8$ the SQP method diverges unless it is warm started.

\begin{table}
  \centering
  \caption{Contact problem 3 (hemisphere), asymmetric manifold:
    solution details of step (0).}
  \label{tab:hemis-asym-s0}
  \sisetup{table-format=1.2e-2,parse-numbers=false}
  \begin{tabular}{S[table-format=1.1]S@{\isep}cSSSS[table-format=2]}
    \toprule
    {$\gamma$} & {NLP} & {active} & {SQP} & {QP} & {constraint} & {SQP} \\
    & {objective} & {set $\setA$} & {tolerance}
    & {tolerance} & {violation} & {iter.} \\
    \midrule
    1   & 6.33e-28 & $\0$          & \E07 & \tQP & 2.95e-08 &  3 \\
    2   & 1.50e-26 & $\0$          & \E07 & \tQP & 1.47e-08 &  4 \\
    2.5 & 2.13e-02 & $\set{13}$    & \E07 & \tQP & 1.49e-08 &  4 \\
    3   & 6.67e-02 & $\set{12,13}$ & \E04 & \E05 & 9.85e-06 &  4 \\
    3.5 & 9.05e-02 & $\set{12,13}$ & \E06 & \E07 & 9.97e-08 &  5 \\
    4   & 6.89e-02 & $\set{13}$    & \E06 & \E07 & 9.88e-08 &  5 \\
    5   & 7.71e-25 & $\0$          & \E07 & \tQP & 1.44e-08 & 11 \\
    6   & 2.72e-01 & $\0$          & \E06 & \E07 & 9.85e-08 & 10 \\
    \bottomrule
  \end{tabular}
\end{table}

\begin{table}
  \centering
  \caption{Contact problem 3 (hemisphere), asymmetric manifold:
    solution details of steps (1) to (3).
    Step (3) on the last line marked with * is warm started from step (1).}
  \label{tab:hemis-asym}
  \sisetup{table-format=1.2e-2,parse-numbers=false}
  \begin{tabular}{S[table-format=1.1]cSSSScS[table-format=2]}
    \toprule
    {$\gamma$}
    & step & {NLP} & {SQP} & {QP} & {constraint} & {compl.} & {SQP} \\
    & & {objective}
    & {tolerance} & {tolerance} & {violation} & {product} & {iter.} \\
    \midrule
    2.5 & (1) & 1.64e-03 & \E07 & \tQP & 1.48e-08 & 0      &  3 \\
        & (2) & 8.62e-33 & \E07 & \tQP & 1.49e-08 & 0.0330 &  3 \\
        & (3) & 1.34e-13 & \E07 & \tQP & 1.47e-08 & 0.0330 &  1 \\
        & (3) & 1.42e-14 & \E10 & \E11 & 9.24e-12 & 0.0330 & 26 \\
    \addlinespace
    3   & (1) & 5.84e-03 & \E06 & \E07 & 4.68e-07 & 0      &  3 \\
        & (2) & 4.96e-27 & \E06 & \E07 & 9.91e-08 & 0.3232 &  3 \\
        & (3) & 3.77e-12 & \E07 & \tQP & 1.49e-08 & 0.3232 &  1 \\
        & (3) & 2.63e-35 & \E10 & \E11 & 9.94e-12 & 0.3232 & 21 \\
    \addlinespace
    3.5 & (1) & 2.33e-02 & \E06 & \E07 & 9.99e-08 & 0      &  4 \\
        & (2) & 3.08e-24 & \E05 & \E06 & 9.96e-07 & 3.4321 &  4 \\
        & (3) & 1.04e-09 & \E06 & \E07 & 9.98e-08 & 3.4320 &  1 \\
        & (3) & 7.66e-27 & \E10 & \E11 & 9.92e-12 & 3.4319 &  3 \\
    \addlinespace
    4   & (1) & 5.81e-02 & \E06 & \E07 & 9.98e-08 & 0      &  5 \\
        & (2) & {--} & {--} & {diverged} & {--} & {--} \\
        & (3) & 5.41e-02 & \E03 & \E04 & 5.10e-05 & 0.0538 &  1\rlap* \\
    \bottomrule
  \end{tabular}
\end{table}

In this example we observe again that the active set of step (0)
changes in subsequent steps for all studied values of $\gamma$.
\Cref{tab:hemis-asym-A} lists the active bounds encountered in each step.

\begin{table}
  \centering
  \caption{Contact problem 3 (hemisphere), asymmetric manifold:
    active sets $\setA$.}
  \label{tab:hemis-asym-A}
  \begin{tabular}{ccc@{\qquad}ccc@{\qquad}ccc@{\qquad}ccc}
    \toprule
    $\gamma$ & step & $\setA$ &
    $\gamma$ & step & $\setA$ &
    $\gamma$ & step & $\setA$ &
    $\gamma$ & step & $\setA$ \\
    \midrule
    $2.5$ & (0) & $\set{13}$ &
    $  3$ & (0) & $\set{12,13}$ &
    $3.5$ & (0) & $\set{12,13}$ &
    $  4$ & (0) & $\set{13}$ \\
         & (1) & $\set{12,13}$ &
         & (1) & $\range{10}{13}$ &
         & (1) & $\range{10}{13}$ &
         & (1) & $\set{12,13}$ \\
         & (2) & $\set{12}$ &
         & (2) & $\set{7}$  &
         & (2) & $\set{7}$  &
         & (2) & -- \\
         & (3) & $\set{12}$ &
         & (3) & $\set{7}$  &
         & (3) & $\set{7}$  &
         & (3) & $\set{12}$  \\
    \bottomrule
  \end{tabular}
\end{table}

\section{Conclusion}
\label{sec:conclusion}
In this article we extended the hybrid DDCM formulation
by adding geometric inequality constraints,
which yields a contact problem in form of a
mathematical program with complementarity constraints.
This MPCC is more complex than the NLP obtained without inequalities
but still substantially less complex than the discrete-continuous NLP
of the direct DDCM approach.
For solving the MPCC we developed a heuristic quick shot approach
that requires at most four NLP solves.
In our extensive numerical experiments
with three different obstacles and two material laws
the quick shot approach is often successful;
in these cases it produces valid solutions of the contact problem.
The majority of unseccussful cases occur for a
non-intuitive asymmetric material law.
In these cases the results could still be used to warm-start
established MPCC methods, which solve a sequence of NLPs
while driving a regularization parameter to zero.
As a side issue, we also discussed robust hybrid DDCM formulations;
these can be used with or without contact
and they remain in the original problem class (NLP or MPCC).
In summary, the article demonstrates that our hybrid DDCM approach
provides a highly flexible universal framework for simulations
in nonlinear elasticity that leads to moderately complex
optimization problems even when addressing contact problems.

\section*{Acknowledgements}

Cristian G.\ Gebhardt gratefully acknowledges the financial support
from the European Research Council through the ERC Consolidator Grant
“DATA-DRIVEN OFFSHORE” (Project ID 101083157).
Marc C.\ Steinbach gratefully acknowledges the financial support
from the Deutsche Forschungsgemeinschaft (DFG, German Research Foundation) --
SFB1463 -- 434502799.

\section*{Appendix: Relations of Strains and Stresses}
Strain and stress tensors appear in physical dimension $d \in \set{1, 2, 3}$.
They are naturally represented by symmetric matrices $E, S \in \Sym{d}$,
which we map to vectors of $n_d$ relevant entries
using the Voigt notation~\cite{Helnwein:2001}.
Here we give the mapping for the strain tensor $E$;
the stress tensor $S$ is treated similarly:
\begin{align*}
  d &= 1\: &
  E &= \mat[1]{e_{11}}, &
  e &= (e_{11}), &
  n_1 &= 1, \\
  d &= 2\: &
  E &= \mat[2]{e_{11} & e_{21} \\ e_{21} & e_{22}}, &
  e &= (e_{11}, e_{22}, e_{21})\tp, &
  n_2 &= 3, \\
  d &= 3\: &
  E &= \mat[3]{
    e_{11} & e_{21} & e_{31} \\
    e_{21} & e_{22} & e_{32} \\
    e_{31} & e_{32} & e_{33}
  }, &
  e &= (e_{11}, e_{22}, e_{33}, e_{32}, e_{31}, e_{21})\tp, &
  n_3 &= 6.
\end{align*}
Now let $\K_d$ denote the set of finite element nodes
on bodies of dimension $d$ and define
$\K \define \K_1 \cup \K_2 \cup \K_3$.
The overall strain and stress vectors $e, s$
are then composed of nodal component vectors $e_k, s_k$, as follows:
\begin{align*}
  e &= (e_k)_{k \in \K} \in \R^{n_e}, &
  s &= (s_k)_{k \in \K} \in \R^{n_s}, &
  n_e = n_s &= \card{\K_1} + 3 \card{\K_2} + 6 \card{\K_3}.
\end{align*}
Next, let $\M_d$ denote the set of different
elastic materials of bodies of dimension $d$.
For each $m \in \M_d$ we then have a constitutive equation given by
\begin{align*}
  \^g_m\: \R^{2 n_d} &\to \R^{n_d}, &
  \^g_m(e, s) &= 0,&
  \rank \^g'(e, s) = n_d.
\end{align*}
The constitutive equation at node $k \in \K_d$
with material $m(k) \in \M_d$ now reads
\begin{align*}
  g_k(e_k, s_k) &= \^g_{m(k)}(e_k, s_k) = 0,
\end{align*}
and the overall map for the entire hybrid system is
\begin{align*}
  g &= (g_k)_{k \in \K}\: \R^{2 n_e} \to \R^{n_e}, &
  g(e, s) &= 0, &
  \rank g'(e, s) &= n_e.
\end{align*}
Considering finally the DDCM context,
we have a distinct data set $\ds_m$ for each $m \in \M_d$.
Corresponding to the structure of $g$,
we then assign to node $k \in \K_d$ with material $m(k) \in \M_d$
the data set $\ds_k \define \ds_{m(k)}$,
and we define the overall data set as the cartesian product
\begin{align*}
  \ds &\define \prod_{k \in \K} \ds_k.
\end{align*}

\bibliographystyle{siam}
\bibliography{ddcm}

\end{document}